\documentclass[12pt]{amsart}
\usepackage{amsfonts,color,latexsym,amssymb,amsmath,epsfig,rotating,array,mathrsfs}

\newtheorem{mahler}{Mahler's Theorem}

\newtheorem{landau}{Landau's Inequality}

\newtheorem{gelfond}{Gelfond's Lemma}

\newtheorem{dfn}{Definition}[section]
\newtheorem{rem}[dfn]{Remark} 
\newtheorem{prop}[dfn]{Proposition}

\newtheorem{thm}[dfn]{Theorem} 
\newtheorem{algor}[dfn]{Algorithm} 
 
\newtheorem{lemma}[dfn]{Lemma}
\newtheorem{cor}[dfn]{Corollary}

\newtheorem{ex}[dfn]{Example}
\definecolor{purple}{rgb}{.5,0,.5}
\definecolor{red}{rgb}{.6,0,0} 
\definecolor{green}{rgb}{0,.5,0} 

\renewcommand{\qed}{$\blacksquare$}

\newcommand{\floor}[1]{\left\lfloor#1\right\rfloor} 
\newcommand{\ceil}[1]{\left\lceil#1\right\rceil}  

\makeatletter
\renewcommand*\env@matrix[1][c]{\hskip -\arraycolsep
  \let\@ifnextchar\new@ifnextchar
  \array{*\c@MaxMatrixCols #1}}
\makeatother

\newcommand{\thth}{^{\text{\underline{th}}}}

\newcommand{\np}{{\mathbf{NP}}}

\newcommand{\pp}{\mathbf{P}}

\newcommand{\eps}{\varepsilon}

\newcommand{\F}{\mathbb{F}}
\newcommand{\Q}{\mathbb{Q}}
\newcommand{\R}{\mathbb{R}}

\newcommand{\C}{\mathbb{C}}

\newcommand{\N}{\mathbb{N}}
\newcommand{\Z}{\mathbb{Z}}
\newcommand{\bO}{\mathbf{O}}

\newcommand{\Zn}{\Z^n}

\newcommand{\Rn}{\R^n}

\newcommand{\Cs}{\C^*}

\newcommand{\Rs}{\R^*}
 
\newcommand{\sA}{\mathscr{A}} 
\newcommand{\cA}{\mathcal{A}} 
\newcommand{\cB}{\mathcal{B}} 
 
\newcommand{\cD}{\mathcal{D}} 
\newcommand{\cE}{\mathcal{E}}

\newcommand{\sL}{\mathscr{L}}

\newcommand{\hA}{\widehat{\cA}} 
\newcommand{\cN}{\mathcal{N}}

\newcommand{\Rsn}{{(\R^*)}^n}
\newcommand{\Csn}{{(\C^*)}^n}

\newcommand{\dia}{$\diamond$}

\newcommand{\supp}{\mathrm{Supp}}

\newcommand{\sign}{\mathrm{sign}}

\author{J.\ Maurice Rojas}
\email{rojas@math.tamu.edu} 
\address{TAMU 3368, College Station, TX \ 77843-3368} 
\thanks{Partially supported by NSF grant CCF-1900881} 
\title[Counting Real Roots in Polynomial-Time for Circuit Systems]
{\mbox{}\\  
\vspace{-1.1in}Counting Real Roots in Polynomial-Time
for Systems Supported on Circuits} 

\begin{document}

\dedicatory{In memory of Tien-Yien Li, June 28, 1945 -- June 25, 2020.} 

\begin{abstract} 
Suppose $A\!=\!\{a_1,\ldots,a_{n+2}\}\!\subset\!\Zn$ has cardinality 
$n+2$, with all the coordinates of the $a_j$ having absolute value 
at most $d$, and the $a_j$ do {\em not} all lie in the same affine hyperplane.  
Suppose $F\!=\!(f_1,\ldots,f_n)$ is an $n\times n$ polynomial system with 
generic integer coefficients at most $H$ in absolute value, and $A$ the union 
of the sets of exponent vectors of the $f_i$. We give the first algorithm 
that, for any {\em fixed} $n$, counts exactly the number of real roots of $F$ 
in time polynomial in $\log(dH)$. 
\end{abstract}

\maketitle 

\vspace{-.8cm}
\section{Introduction} 
Solving sparse polynomial systems remains a 
challenging problem, even 40 years after the dawn of fewnomial theory 
\cite{kho,few}. More recently, connections have emerged 
between fewnomial theory over finite fields, cryptography, 
and number theory \cite{cfklls,bomb,cgrw}; and 
sparse polynomial systems over the real numbers are important in applications 
including computational biology and biochemistry \cite{sparsebio1,sparsebio2,
sparsechem1, sparsechem2}, and circuit complexity \cite{realtau}. However, 
efficiently counting the number of real 
roots, or even just finding a reasonably tight upper bound on the number of 
real roots, remains an open problem. Here, we focus on the problem of exactly 
counting real roots, and roots in any given orthant.
In what follows, all $O$-constants are absolute and effective, 
{\em time} will refer to the number of (deterministic) bit operations in  
the classical Turing model of computation, and we will use $\#S$ for the 
cardinality of a set $S$. 
 
Suppose $A\!=\!\{a_1,\ldots,a_t\}\!\subset\!\Zn$ has cardinality $t$, 
with all coordinates of the $a_j$ having absolute value at most $d$. 
Writing $x^{a_j}\!:=\!x^{a_{1,j}}_1 \cdots x^{a_{n,j}}_n$ and 
$f(x)\!=\!\sum^t_{j=1}c_j x^{a_j}\!\in\!\Z\!\left[x^{\pm 1}_1,\ldots,
x^{\pm 1}_n\right]$, we define the {\em support} of 
$f$ to be $\supp(f)\!:=\!\{a_j \; | \; c_j\!\neq\!0\}$.  
We then call a system of the form $F\!:=\!(f_1,\ldots,f_n)
\!\in\!\Z\!\left[x^{\pm 1}_1,\ldots,x^{\pm 1}_n\right]^n$,   
with $f_i(x)\!:=\!\sum^t_{j=1} c_{i,j}x^{a_j}$ for all $i$ and 
$\#\bigcup^n_{i=1} \supp(f_i)\!=\!t$, a {\em $t$-nomial $n\times n$ 
system (over $\Z$) supported on $A$}. We denote the positive orthant by 
$\Rn_+$ and call a root of $F$ in $\Rn_+$ a {\em positive} root. We also let 
$\Rs\!:=\!\R\setminus\{0\}$. 

If the $a_j$ do not all lie in the same affine hyperplane then we clearly have 
$t\!\geq\!n+1$. It is\linebreak 
\scalebox{.97}[1]{natural to assume that the exponent vectors are non-coplanar 
in this sense, and we will do so,}\linebreak 
\scalebox{.985}[1]{for otherwise one could use a monomial change of 
variables to reduce $F$ to a system in fewer}\linebreak 
variables: See Remark \ref{rem:mono} from Section 
\ref{sec:back} below. Our main theorem gives a dramatic speed-up for counting 
the exact number of real roots of $F$ in the special case $t\!\leq\!n+2$. 
\begin{thm}
\label{thm:main}  
Following the notation above, assume the coefficient matrix 
$[c_{i,j}]$ of $F$ is generic and $t\!\leq\!n+2$. Then, in time 
$O\!\left(n\log(nd)+n^2\log(nH)\right)^{2n+4}$,  
we can determine the number of roots of $F$ in $\Rn$, $\Rsn$, and $\Rn_+$.  
Furthermore, if $t\!=\!n+1$, then we can do the same in time 
$O\!\left(n^{3.373}\log^2(nH)\right)$.  
\end{thm} 

\noindent 
We prove Theorem \ref{thm:main} in Section \ref{sub:proof}, based 
mainly on Algorithms \ref{algor:sign} and \ref{algor:main} from Section 
\ref{sec:algor}. A central ingredient is diophantine approximation.  

The genericity condition for counting roots in $\Rsn$ and $\Rn_+$, when 
$t\!=\!n+2$, is the non-singularity of a collection of sub-matrices of 
$[c_{i,j}]$, depending solely on the support of $F$, and is checkable in time 
$O\!\left(n^{3.373}\log^2(ndH)\right)$: See Lemma \ref{lemma:gale} 
and Corollary \ref{cor:gen} of Section \ref{sub:n+2}. In 
particular, the fraction of coefficient matrices failing to satisfy this 
genericity condition is no greater than $\frac{2n^2}{2H+1}$. The genericity 
condition for counting affine roots is more technical but still holds 
practically often: probability $1-\eps$ when $H$ has $\Omega(n\log(d)
-\log\eps)$ bits (see Section \ref{sub:proof}). Root counting without 
genericity assumptions is rather non-trivial: Deciding infinitude for the 
number of real (or positive) roots in 
time polynomial in $(n\log(dH))^n$, when $t\!=\!n+2$ and 
$f_2\!=\cdots\!=\!f_n$, is already an open question \cite{brs,bihan}. 
Furthermore, for any fixed $\eps\!>\!0$, deciding whether $F\!=\!(f_1,
\ldots,f_1)$ has any real (or positive) 
roots is $\np$-hard already for $t\!=\!n+n^\eps$ \cite{brs}.
\begin{rem} 
We count roots without multiplicity. In particular, degenerate\footnote{roots 
yielding a Jacobian with less than full rank.} isolated 
roots are not a problem, and are counted correctly by our algorithms. \dia 
\end{rem} 

Other than an algorithm for the very special case $(n,t)\!=\!(1,3)$ 
from \cite{brs}, the best previous complexity bound for $t\!=\!n+2$ 
appears to have been $n^{(1+n/2)^2}d^{O(n^2)}(\log H)^{O(n^2)}$,  
as a consequence of more general algorithms (see, e.g., \cite{bprdim,bpr}), 
based on older computational algebra techniques 
(see, e.g., \cite{chigo,bkr,canny,renegar}). (One can also derive a 
$(d\log H)^{O(n)}$ bound via \cite{proy} if one assumes the complex 
zero set is finite.) There have also been important recent advances from 
the point of view of numerical conditioning (e.g., \cite{krick1,krick2}), 
even enabling practical computation of homology of real projective sets,  
but work in this direction has not yet focussed on speed-ups like Theorem 
\ref{thm:main}. With few exceptions, earlier work on solving polynomial 
systems over the real numbers focused on coarser complexity bounds that 
ignored the finer monomial term structure. 
\begin{ex} 
\label{ex:rigged} 
Consider the $7$-nomial $5\times 5$ system $F\!=\!(f_1,\ldots,f_5)$ defined 
by\\ 
\scalebox{.73}[.73]{\vbox{
\begin{eqnarray*} 
\left(2x^{36}_1 x^{194}_2 x^{50}_3 x^{82}_4 x^{60}_5 
+x^{76}_1 x^{240}_2 x^{41}_4 x_5 
+x^{74}_1 x^{179}_2 x^{25}_3 x^{57}_5 
+x^{25}_1 x^{203}_2 x^{44}_3 x_4  
+x^{20}_1 x^{167}_2 x^{64}_3 x^{12}_4 x^{68}_5 
-37137cx^{58}_1 x^{194}_2 x^{24}_3 x^{36}_4 x^{25}_5 
-\text{\scalebox{.8}[.8]{$\frac{9}{2}$}}x^{166}_3 x^{68}_4 x^{343}_5,\right. \\ 
x^{36}_1 x^{194}_2 x^{50}_3 x^{82}_4 x^{60}_5 
+2x^{76}_1 x^{240}_2 x^{41}_4 x_5 
+x^{74}_1 x^{179}_2 x^{25}_3 x^{57}_5 
+x^{25}_1 x^{203}_2 x^{44}_3 x_4  
+x^{20}_1 x^{167}_2 x^{64}_3 x^{12}_4 x^{68}_5 
-24849cx^{58}_1 x^{194}_2 x^{24}_3 x^{36}_4 x^{25}_5 
-\text{\scalebox{.8}[.8]{$\frac{21}{4}$}}x^{166}_3 x^{68}_4 x^{343}_5,\\  
x^{36}_1 x^{194}_2 x^{50}_3 x^{82}_4 x^{60}_5 
+x^{76}_1 x^{240}_2 x^{41}_4 x_5 
+2x^{74}_1 x^{179}_2 x^{25}_3 x^{57}_5 
+x^{25}_1 x^{203}_2 x^{44}_3 x_4  
+x^{20}_1 x^{167}_2 x^{64}_3 x^{12}_4 x^{68}_5 
-21009cx^{58}_1 x^{194}_2 x^{24}_3 x^{36}_4 x^{25}_5 
-\text{\scalebox{.8}[.8]{$\frac{21}{4}$}}x^{166}_3 x^{68}_4 x^{343}_5,\\  
x^{36}_1 x^{194}_2 x^{50}_3 x^{82}_4 x^{60}_5 
+x^{76}_1 x^{240}_2 x^{41}_4 x_5 
+x^{74}_1 x^{179}_2 x^{25}_3 x^{57}_5 
+2x^{25}_1 x^{203}_2 x^{44}_3 x_4  
+x^{20}_1 x^{167}_2 x^{64}_3 x^{12}_4 x^{68}_5 
-20769cx^{58}_1 x^{194}_2 x^{24}_3 x^{36}_4 x^{25}_5 
-\text{\scalebox{.8}[.8]{$\frac{21}{4}$}}x^{166}_3 x^{68}_4 x^{343}_5,\\  
\left. x^{36}_1 x^{194}_2 x^{50}_3 x^{82}_4 x^{60}_5 
+x^{76}_1 x^{240}_2 x^{41}_4 x_5 
+x^{74}_1 x^{179}_2 x^{25}_3 x^{57}_5 
+x^{25}_1 x^{203}_2 x^{44}_3 x_4  
+2x^{20}_1 x^{167}_2 x^{64}_3 x^{12}_4 x^{68}_5 
-20754cx^{58}_1 x^{194}_2 x^{24}_3 x^{36}_4 x^{25}_5 
-\text{\scalebox{.8}[.8]{$\frac{21}{4}$}}x^{166}_3 x^{68}_4 x^{343}_5\right).   
\end{eqnarray*}}}\\  
Then Algorithm \ref{algor:main} from Section \ref{sec:algor} (simulated 
in a few lines of {\tt Maple} code\footnote{Using {\tt Maple 2019} on a 
Dell XPS 13 laptop with an Intel core i7-5500u microprocessor, 8 Gb RAM, and 
a 256Gb solid state hard-drive, running Ubunu 19.10. {\tt Maple} 
code available on request.}) tells us in under a second that 
$F$ has exactly $2$, $6$, $6$, $2$, $2$, or $0$ positive roots, respectively 
when $c$ is $\frac{1}{20731}$, $\frac{1}{20730}$, $\frac{1}{14392}$, 
$\frac{1}{14391}$, $\frac{1}{13059}$, or $\frac{1}{13058}$. 
(All roots in $(\Rs)^5$ of these $F$ happen to lie in $\R^5_+$.) 
We will return to this family in Section \ref{sub:n+2}, and see another 
example there as well. It is interesting to observe that 
{\tt Maple}'s {\tt Solve} command (which employs Gr\"{o}bner 
bases) gives no useful information about any of these systems, even after 
$3$ hours.$^2$ {\tt Bertini} (a state of the art homotopy solver, version 
1.4 \cite{bertini}), on each of the preceding systems, immediately returns a 
message stating  ``ERROR: The system is 
numerically zero 0! Please input a non-degenerate system. Bertini will now 
exit due to this error.'' This is partially because each of our $F$ 
has\footnote{ Via Kushnirenko's Theorem \cite{kush,why}, 
Ioannis Emiris' {\tt MixVol} code \cite{emiris}, and a simple 
check that the underlying facial systems have no roots in $(\Cs)^5$.} 
over $245$ million 
roots in $(\Cs)^5$,  
and older polynomial system solving techniques have complexity super-linear,  
or worse, in the number of complex roots. \dia 
\end{ex}  

\begin{rem} 
The main intent of Theorem \ref{thm:main} is to set the 
stage (building on the framework of \cite{ppr,epr1,epr2}) for more practical 
improvements in real-solving, such as complexity 
sub-exponential in $n$ in the average-case/smoothed analysis setting for 
sparse systems. In particular, just as binomial systems are a building block 
for polyhedral homotopy algorithms for arbitrary $n\times n$ systems 
\cite{hs,verschelde,leeli}, $(n+2)$-nomial $n\times n$ systems are a building 
block behind recent optimization techniques such as SONC/SAGE-optimization 
(see, e.g., \cite{snc,chandra,timo}). So while tackling the remaining 
exceptional cases (involving infinitely many real roots in $\Rsn$) is 
important, such cases are provably rare when the coefficients are random. \dia 
\end{rem} 

There has been growing interest in generalizing {\em Descartes' Rule of Signs} 
(see, e.g., \cite{descartes,grabiner}) from univariate polynomials to 
$n\times n$ polynomial systems. This began with\linebreak  
Khovanski's seminal {\em Theorem on Real Fewnomials} \cite{few} which, 
in our notation, asserted an upper bound of 
$\displaystyle{2^{\binom{t}{2}}(n+1)^t}$ for the number of non-degenerate 
positive roots of any $t$-nomial $n\times n$ system. It was then shown 
in \cite{tri} that Khovanski's bounds could be greatly 
improved for various structured systems, e.g., the correct tight upper bound 
on the number of isolated\footnote{even allowing degenerate isolated roots} 
positive roots for $2\times 2$ systems of trinomials is $5$ --- far less 
than the best previous bound of $248832$. Sharper upper bounds for 
new families of systems, including a tight upper bound of $n+1$ 
(resp.\ $(n+1)2^n$) roots in $\Rn_+$ (resp.\ $\Rsn$) 
for the case $t\!=\!n+2$, were then derived 
in \cite{bbs}. Explicit families of systems 
attaining these bounds for each $n$ were then given in \cite{pr} (see 
also \cite{bihan}). 
Khovanski's general upper bound was vastly improved to 
$\displaystyle{\frac{e^2+3}{4}2^{\binom{t-n-1}{2}}n^{t-n-1}}$ positive roots 
in \cite{bs}, and a remarkable bound for curve intersections was derived later 
in \cite{lostfew}. More recently, a beautiful and considerably sharper 
{\em average-case} upper bound of\\ 
\mbox{}\hfill $\displaystyle{\frac{1}{2^{n-1}}\cdot \frac{t!}{n!(t-n)!}}$\hfill
\mbox{}\\  
for the number of positive roots was proved in \cite{bet}, using independent 
real Gaussians for the coefficients. It should also be noted that counting 
roots on coordinate subspaces quickly abuts $\#\pp$-hardness, already 
for $n\times n$ binomial systems: See \cite{catdick,mondal} and Remark 
\ref{rem:sharpp} in Section \ref{sub:bino} below. So it wasn't just the 
convenience of torus actions that made earlier work on fewnomials focus on 
$\Rsn$ instead of $\Rn$. 

However, fewnomial bounds so far have not made significant use of the signs of 
the coefficients (much less their values) and such refined bounds remain 
elusive: See, e.g., \cite[Thm.\ 2.1]{batesgale} and \cite{bihandickenstein,
signone,optcktdescartes}. The latter works, particularly 
\cite{optcktdescartes}, 
culminated in a refined characterization of the maximal number of positive 
roots --- incorporating the signs of $n\times n$ sub-determinants of the 
coefficient matrix $[c_{i,j}]$ and the matroidal structure of $\cA$ --- in the 
case $t\!=\!n+2$. Nevertheless, no algorithm for {\em exactly} counting the 
real or positive roots, faster than combining more general results on rational 
univariate reduction (see, e.g., \cite{kronecker,tgcp,rouillier}) with the 
computation of real dimension (see, e.g., \cite{bprdim}) or real root isolation 
(see, e.g., \cite{sagraloff}), appears to have been known before 
Theorem \ref{thm:main} above. 

Exactly counting the real or positive roots of $F$, and even formulating 
a reasonable generalization of Descartes' Rule, appears to be much 
harder for $t\!\geq\!n+3$. This is why there is so much attention on the case 
$t\!=\!n+2$ to develop further intuition. An even harder open question is 
the complexity of actually {\em approximating} the real roots of such $F$
and we hope to address this in future work.  

Our main tools are reduction to a canonical form (a special case of 
{\em Gale Dual Form} from \cite{bs}) and a careful application of   
diophantine approximation to the critical values  
of this reduction. In particular, the locus of $F$ with 
degenerate roots forms a {\em discriminant variety} which partitions the 
coefficient space into connected open regions 
we call {\em chambers} (see, e.g., \cite[Ch.\ 11]{gkz94}). 
Classical topological results, such as {\em Hardt's Triviality Theorem} 
\cite{hardt}, tell us that counting the 
real roots of $F$ is tantamount to identifying the chamber 
in which $F$ lies. Such a calculation is challenging, since the theory of 
$\cA$-discriminants \cite{gkz94} does not directly provide a tractable 
description of our chambers. However, a combination of Rolle's Theorem 
with Gale Dual Form allows one to replace chamber identification by the 
determination of signs of the critical values and poles of a single univariate 
rational function. 

A new obstacle is that the usual univariate root-finding algorithms, 
combined with classical height bounds on polynomial roots,  
do not yield a useful complexity bound. Indeed, the degree of the resulting 
univariate reduction can be so high that a naive use of real root isolation 
would lead to complexity super-linear in $n^{n/2} d^n$. So we leverage the 
special structure of the {\em derivative} of our univariate reduction to apply 
a powerful theorem from diophantine approximation: A refinement of an estimate 
of Baker and Wustholtz on linear forms in logarithms of algebraic numbers (see 
\cite{baker,bakerwustholtz,matveev,bms} or Section \ref{sec:back} below).

\begin{rem} 
\label{rem:abc} 
A by-product of our framework is that sufficiently 
sharp lower bounds for linear forms in logarithms of algebraic numbers 
would easily imply that root counting in $\Rn$, for generic systems over 
$\Z$ supported on {\em circuits}\footnote{...not to be confused with the 
circuits from complexity theory (which are layered directed graphs with 
specially labelled nodes having additional structure). See Section 
\ref{sub:n+2}.}, can be sped up to worst-case complexity {\em polynomial} in 
$n$. While the necessary new diophantine bounds 
appear far out of reach, Baker has proved \cite{bakerabc}, in the special case 
of linear forms of logarithms of rational numbers, 
that such bounds would follow from a refined version of the famous 
Masser-Oesterle $abc$-Conjecture \cite{masser,oesterle,nitaj}. Unfortunately, 
the latter refinement also 
appears out of reach (as of early 2021). This is one reason that new 
average-case speed-ups, using geometric numerical conditioning techniques 
(e.g., \cite{epr1,epr2}) instead of diophantine approximation, may arrive 
sooner than new worst-case speed-ups. \dia 
\end{rem} 

\section{Background} 
\label{sec:back}  
\subsection{The Complexity of Linear Algebra over $\Z$} 
\label{sub:lin} 
Let $\omega$ denote the well-known matrix multiplication exponent, i.e., 
the infimum over all $\omega$ such that there exists an algorithm 
that can multiply an arbitrary pair of $n\times n$ matrices, in any field $K$, 
using $O(n^\omega)$ field operations in $K$. The best current upper bound 
is $\omega\!<\!2.3728639$ \cite{virgi,legall}.  Recall the notions of 
{\em reduced row echelon form} and {\em leading entries} of a matrix, from 
basic linear algebra 
(see, e.g., \cite{prasolov}). For any rational number $\frac{p}{q}$ with 
$p,q\!\in\!\Z$ and $\gcd(p,q)\!=\!1$, its {\em (absolute) logarithmic height} 
is $h(p/q)\!:=\!\max\{\log |p|,\log |q|\}$. (We set $h(0)\!:=\!0$.) We will 
first need a result on the bit complexity of row reduction for matrices: 
\begin{dfn} \cite{hermite,schrijver} 
We call a matrix $U\!\in\!\Z^{n\times n}$ with determinant $\pm 1$
{\em unimodular}. Given any matrix $M\!\in\!\Z^{n\times t}$, we then call
any identity of the form $UM\!=\!R$, with $U\!\in\!\Z^{n\times n}$ unimodular
and $R$ upper-triangular with all leading
entries positive, a {\em Hermite factorization}. Finally, we call
any identity of the form $UMV\!=\!S$, with $U\!\in\!\Z^{n\times n}$
and $V\!\in\!\Z^{t\times t}$ both unimodular, and $S$ with diagonal
entries $s_1,\ldots,s_n$ satisfying $s_1|s_2, \ldots,  
s_{n-1}|s_n$, a {\em Smith factorization of $M$}. \dia
\end{dfn}
\begin{lemma} 
\label{lemma:lin} 
Suppose $M\!\in\!\Z^{n\times t}$ with $t\!\geq\!n$ and all entries having 
absolute value at most $H$. Then, in time $O\!\left(tn^{\omega}\log^2(nH)
\right)$, we can row-reduce $M$ to reduced row echelon form 
$R\!\in\!\Q^{n\times t}$ with every nonzero entry of $R$ having logarithmic 
height $O(n\log(nH))$. \qed 
\end{lemma} 

\noindent 
The complexity bound from Lemma \ref{lemma:lin}
is easily obtainable by applying a fast algorithm for Hermite Factorization 
(see, e.g., \cite[Ch.\ 6]{storjophd}) to reduce $M$ to an upper-triangular 
matrix in $\Z^{n\times t}$, dividing through by the leading entries, and then
back-substituting to obtain the desired reduced row echelon form.
An illuminating alternative discussion, via the {\em Newton}\linebreak 
\scalebox{.91}[1]{{\em identities} for
sums of powers of roots of a polynomial, can be found in \cite[Ch.\ 15,
Sec. 15.5]{bcss}.} 

Via Cramer's Rule and Hadamard's classical inequality on the absolute values 
of determinants \cite[Thm.\ 1]{mignotte}, we can easily obtain the 
following related bound:
\begin{lemma}
\label{lemma:rightnull}
If $\cA\!\in\!\Z^{n\times (n+1)}$ has rank $n$ and the entries of the $i$th row
of $\cA$ have absolute value at most $d_i$, then any generator 
$(b_1,\ldots,b_{n+1})^\top\!\in\!\Z^{(n+1)\times 1}$ 
of the right-null space of $\cA$, with $\gcd(b_1,\ldots,b_{n+1})\!=\!1$,
satisfies $|b_j|\!\leq\!n^{n/2}\prod\limits^n_{i=1} d_i$ for all $j$. \qed
\end{lemma}

We will also need the following complexity bound on Smith factorization: 
\begin{thm} \label{thm:smith} \cite[Ch.\ 8]{storjophd} 
Suppose $M\!\in\!\Z^{n\times k}$ (with $k\!\leq\!n$) 
has all entries of absolute value 
at most $d$. Then a Smith factorization $UMV\!=\!S$ for $M$ can be found in 
time $O\!\left(kn^{\omega}\log^2(nd)\right)$, with all the entries 
of $U,V,S$ having logarithmic height $O(n\log(nd))$. \qed 
\end{thm} 

\subsection{Binomial and $(n+1)$-nomial Systems over $\Rsn$} 
\label{sub:bino} 
\scalebox{.9}[1]{A simple, folkloric algebraic/analytic}\linebreak 
fact we'll need is the following:  
\begin{prop}
\label{prop:mono} 
Suppose $\cA,\cB\!\in\!\Z^{n\times n}$ and $x\!=\!(x_1,\ldots,x_n)$ is a 
vector of indeterminates. Let us define $x^{\cA}$ to be the 
vector of monomials $\left(x^{a_{1,1}}_1\cdots x^{a_{n,1}}_n,
\ldots,x^{a_{1,n}}_1\cdots x^{a_{n,n}}_n\right)$, where 
$\cA\!=\![a_{i,j}]$. Then 
$(x^{\cA})^{\cB}\!=\!x^{\cA\cB}$ and, if $\cA$ is 
unimodular, the function defined by $x\mapsto x^{\cA}$ defines\linebreak 
\scalebox{.9}[1]{an analytic 
group automorphism of $\Csn$ that restricts to an analytic group automorphism 
of $\Rn_+$. \qed}  
\end{prop} 
\begin{rem}
\label{rem:mono} 
A simple consequence of Proposition \ref{prop:mono} is that 
if $f\!\in\!\R\!\left[x^{\pm 1}_1,\ldots,x^{\pm 1}_n\right]$ is an 
$n$-variate $t$-nomial with support $A$, and $d$ is the dimension 
of the smallest affine subspace containing $A$, then there is a {\em 
monomial change of variables} $x\!=\!y^U$ (with $U$ unimodular), and a 
monomial $y^b\!\in\!\R\!\left[y^{\pm 1}_1,\ldots,y^{\pm 1}_d\right]$, 
such that $g(y)\!:=\!y^b f\!\left(y^U\right)\!\in\!\R\!\left[y^{\pm 1}_1,
\ldots,y^{\pm 1}_d\right]$ is a $d$-variate $t$-nomial, and the zero 
set of $f$ in $\Rsn$ is analytically isomorphic to the Cartesian 
product\linebreak 
\scalebox{.95}[1]{of $(\Rs)^{n-d}$ and the zero set of $g$ in $(\Rs)^d$. So 
$A$ in an affine hyperplane implies that the zero set}\linebreak 
\scalebox{.94}[1]{of $f$ in $\Rsn$ can be easily characterized by the zero set 
of another $t$-nomial in fewer variables. \dia}   
\end{rem}     

Another consequence of Lemma \ref{lemma:lin} is that we can almost trivially 
count the positive roots of {\em binomial} systems, provided 
the exponent vectors are in general position.
\begin{prop}
\label{prop:bino} 
Suppose $c\!=\!(c_1,\ldots,c_n)\!\in\!\Rsn$, 
$a_1,\ldots,a_n\!\in\!\Zn$, $\cA$ is the $n\times n$ matrix\linebreak 
with $j$th column $a_j$ for all $j$, $U\cA V\!=\!S$ is a Smith factorization 
of $\cA$, and $c'\!:=\!\left(c'_1,\ldots,c'_n\right)\!:=\!c^V$. 
Let $s_j$ be the $(j,j)$ entry of $S$. 
Then $G\!:=\!(x^{a_1}-c_1,\ldots,x^{a_n}-c_n)$ 
and $\left(y^{s_1}_1-c'_1,\ldots,y^{s_n}_n-c'_n\right)$ have the same number 
of roots in $\Rn_+$ (resp.\ $\Rsn$, $\Csn$).  
In particular, $G$ has exactly $0$, $1$, or infinitely many roots in 
$\Rn_+$ under the following respective conditions: \\ 
\mbox{}\hspace{1cm}\fbox{$0$}: Some $c_i$ is negative or 
[$\mathrm{Rank}(\cA)\!=\!j\!<\!n$ and $c'_i\!\neq\!1$ for some 
$i\!\in\!\{j+1,\ldots,n\}$]. \\ 
\mbox{}\hspace{1cm}\fbox{$1$}: $c\!\in\!\Rn_+$ and $\det \cA\!\neq\!0$. \\ 
\mbox{}\hspace{.8cm}\fbox{$\infty$}: $c\!\in\!\Rn_+$, 
$\mathrm{Rank}(\cA)\!=\!j\!<\!n$, and $c'_{j+1}\!=\cdots=\!c'_n\!=\!1$. \qed 
\end{prop} 

\noindent 
Proposition \ref{prop:bino} follows directly from Proposition \ref{prop:mono}. 
Both facts are folkloric in the toric geometry/Lie group literature (see, 
e.g., \cite{hs} or \cite[Ch.\ 3]{bombgu}). 
A more in-depth discussion of binomial systems can be 
found in \cite{catdick,tianran,ppr}. 

Counting roots in $\Rsn$ is slightly more complicated but still admits 
efficient formulae. 
\begin{prop} 
\label{prop:rsn} 
Following the notation of Proposition \ref{prop:bino}, assume the 
exponent vectors $a_1,\ldots,a_n$ are linearly independent. Let 
$r$ denote the rank, over the field $\F_2$, of the mod $2$ reduction 
of $\cA$. Also let $\overline{V}$ denote the upper $r\times n$ sub-matrix of 
$V$. Then the map $m : \Rsn \longrightarrow \Rsn$ defined 
by $m(x)\!:=\!x^{\cA}$ is $2^{n-r}$-to-$1$, and the $i$th 
coordinate of its range is $\Rs$ (resp.\ $\R_+$) if and only 
if some (resp.\ no) entry of the $i$th column of $\overline{V}$ is odd. 
In particular, $F$ has exactly $0$ (resp.\ $2^{n-r}$) roots in $\Rsn$ 
if and only if $c'_i\!<\!0$ for some (resp.\ no) $i\!\geq\!r+1$.  
\end{prop}

\noindent 
{\bf Proof:} First note that, by the definition of Smith factorization, 
we have that the diagonal entries $s_i$ of $S$ are such that 
$s_1,\ldots,s_r$ are odd and $s_{r+1},\ldots,s_n$ are even. By 
Proposition \ref{prop:mono}, exponentiating by $U$ or $V$ induces a 
permutation of the open orthants of $\Rn$. In particular, we see that 
the range of $x\mapsto x^S$ is exactly $(\Rs)^r\times \R^{n-r}_+$.  
So the pre-image assertion on $m$ is proved.

Now note that the range of $m$ must be $\left((\Rs)^r\times \R^{n-r}_+\right)
^V$. Since the sign of $\eps^{v_{i,j}}_j$ depends only on the parity of 
$v_{i,j}$ and the sign of $\eps_j$, we then see that the $i$th coordinate 
of the range of $m$ is $\Rs$ (resp.\ $\R_+$) if and only if 
some (resp.\ no) entry of $\{v_{1,i},\ldots,v_{r,i}\}$ is odd. So now 
we know the range of $m$. 

The final remaining assertion follows from our earlier definition 
$c'\!:=\!c^V$ and our earlier assumption that $c\!\in\!\Rsn$. \qed 

\begin{rem} 
\label{rem:sharpp} 
It is important to remember that counting {\em affine} roots of 
$n\times n$ binomial systems, with {\em arbitrary}   
exponents and $n$ varying, is in fact $\#\pp$-complete \cite{catdick}. This 
implies, for instance, that polynomial-time root counting in $\Rn$ for 
systems of the form $(x^{\alpha_1}f_1,\ldots,x^{\alpha_n}f_n)$ ---  
with $F\!=\!(f_1,\ldots,f_n)$ an $(n+1)$-nomial $n\times n$ system and 
$\alpha_1, \ldots,\alpha_n\!\in\!(\N\cup\{0\})^n$ --- would imply 
$\pp\!=\!\np$. In particular, our restriction to $n\times n$ 
systems with {\em union} of supports having cardinality $\leq\!n+2$ is 
critical if we are to have any hope of someday counting affine roots in 
$\pp$. \dia 
\end{rem} 

We can now state more explicitly how we deal with positive root counting 
for $t$-nomial systems in the case $t\!=\!n+1$. 
\begin{lemma}
\label{lemma:n+1} 
If $F\!=\!(f_1,\ldots,f_n)\!\in\!\Z\!\left[x^{\pm 1}_1,\ldots,x^{\pm 1}_n
\right]^n$ is an $(n+1)$-nomial $n\times n$ system, with union of supports 
$A\!=\!\{a_1,\ldots,a_{n+1}\}$ {\em not} lying in an affine hyperplane, and 
the coefficient matrix of $F$ has rank $n$, then the number of positive roots 
of $F$ is either $0$ or $1$. Furthermore, if all the coefficients of all the  
$f_i$ have absolute value at most $H$, then we can determine the number of 
positive roots of $F$ in time $O(n^{1+\omega}\log^2(nH))$.  
\end{lemma} 
\begin{rem} 
\label{rem:disturb} 
The reader disturbed by the complexity bound being independent of $A$ 
may be reassured to know that (a) checking the hyperplane condition 
takes time dependent on $A$ and (b) the analogue 
of our lemma for counting roots in $\Rsn$ (Corollary \ref{cor:n+1+} below) 
has complexity dependent on $A$. \dia 
\end{rem}  

\noindent 
{\bf Proof of Lemma \ref{lemma:n+1}:} By our assumption on the coefficient 
matrix, we may reorder monomials so that the left-most $n\times n$ 
minor of the coefficient matrix has nonzero determinant. So we may divide every 
$f_i$ by $x^{a_{n+1}}$ without changing the roots of $F$ in $\Csn$, 
and assume $a_1,\ldots,a_n$ are linearly independent and $a_{n+1}\!=\!\bO$. 

From Lemma \ref{lemma:lin} it is then 
clear that we can reduce the coefficient matrix of $F$, 
$[c_{i,j}]\!\in\!\Z^{n\times (n+1)}$, to 
a reduced row echelon form in $\Q^{n\times (n+1)}$, in 
time $O(n^{1+\omega}\log^2(nH))$. The underlying linear combinations of 
rows can then be applied to the equations $f_i\!=\!0$ so that $F\!=\!\bO$ can 
be reduced to a binomial system of the form $x^{\cA}\!=\!\gamma$ where 
$\gamma\!=\!(\gamma_1,\ldots,\gamma_n)$, $\cA\!\in\!\Z^{n\times n}$,  
and the solutions of $x^{\cA}\!=\!\gamma$ in $\Csn$ are the same as the roots 
of $F$ in $\Csn$.  

Clearly then, $\gamma_i\!\leq\!0$ for any $i$ implies that $F$ 
has no positive roots. In which case, we simply report that $F$ has 
$0$ positive roots and conclude, having taken time $O(n^{1+\omega}\log^2(nH))$. 

\scalebox{.96}[1]{So $\gamma\!\in\!\Rn_+$ implies that $F$ has exactly 
$1$ positive root by Proposition \ref{prop:bino}, and we are done. \qed} 

\medskip 
A simple consequence of our development so far is a method to efficiently 
count roots in $(\Rs)^n$ for generic $(n+1)$-nomial systems. 
\begin{cor}
\label{cor:n+1+} 
Following the notation and assumptions of Lemma \ref{lemma:n+1} and 
its proof, the number of roots of $F$ in $\Rsn$ is either $0$ or $2^{n-r}$, 
where $r$ is the rank, over the field $\F_2$, of the mod $2$ reduction of 
$\cA$. In particular, we can determine the number of roots of $F$ in $\Rsn$ 
in time $O(n^{1+\omega}\log^2(ndH))$, where $d$ is the maximum absolute 
value of any entry of $A$.  
\end{cor}

\noindent 
{\bf Proof:} Continuing from the proof of Lemma \ref{lemma:n+1} (and 
having already reduced our input $(n+1)$-nomial $n\times n$ system 
to a binomial system), it is clear that Proposition \ref{prop:rsn} tells us 
that we can easily count the roots of $F$ in $\Rsn$: We merely need to check 
the signs of $\gamma'_{r+1},\ldots,\gamma'_n$ where $\gamma'\!:=\!\gamma^V$ 
and $U\cA V\!=\!S$ is a Smith factorization of $\cA$. 
More precisely, instead of 
computing $\gamma^V$, we compute $\sign(\gamma)^{(V \text{ mod } 2)}$. 
Computing the mod $2$ reduction of $V$ takes time $O(n^2)$ and 
then computing the resulting vector of signs clearly takes 
time just $O(n^2)$. So the only remaining work (after performing 
Gauss-Jordan elimination on the coefficient matrix of $F$) is 
extracting the Smith factorization of $\cA$ via, say, Theorem \ref{thm:smith}. 
So our final complexity bound is $O(n^{1+\omega}\log^2(nH)
+n^{1+\omega}\log(dH))$, 
which is clearly bounded from above by our stated bound. \qed

\subsection{Circuits, $(n+2)$-nomial systems, and Gale Dual Form with Heights} 
\label{sub:n+2} 
We now show how to reduce root counting in $\Rsn$ for $F$ to root counting 
in certain sub-intervals of $\R$ for a linear combination of logarithms in one 
variable. This reduction dates back to \cite{bs}, if not earlier, but our 
statement here includes height and computational complexity bounds that appear 
to be new. Before proving our reduction, however, let us recall the 
combinatorial/geometric notion of a {\em circuit}: 
\begin{dfn} We call any subset $A\!=\!\{a_1,\ldots,a_{m+2}\}\!\subset\!\Rn$
with $\#A\!=\!m+2$ a {\em circuit} if and only if the $(n+1)\times (m+2)$ 
matrix $\hA$ with $j$th column $\begin{bmatrix}1\\ a_j\end{bmatrix}$
has right nullspace of dimension one. In which case, we call any
generator $b\!\in\!\Z^{(m+2)\times 1}\setminus\{\bO\}$ for the right 
nullspace of $\hA$, with $1$ for its gcd
of coordinates, a {\em (minimal) circuit relation} of $A$. We also
call $A$ a {\em degenerate} circuit if and only if $b$ has at least one zero
coordinate. \dia
\end{dfn}

\noindent 
Note that all circuit relations for a fixed circuit 
(other than the trivial relation $\bO$) have zero entries 
occuring at the {\em same} set of coordinates. More precisely, the
following proposition is elementary.
\begin{prop}
\label{prop:ckt}
Any circuit $A\!=\!\{a_1,\ldots,a_{m+2}\}\!\subset\!\Zn$ has a
{\em unique} subset $\Sigma\!=\!\{a_{i_1},\ldots,a_{i_{\ell+2}}\}$ with 
$\Sigma$ a non-degenerate circuit of cardinality $\ell+2$. In particular, 
$\{i_1,\ldots,i_{\ell+2}\}$
is exactly the set of indices of the nonzero coordinates of any (non-trivial) 
circuit relation for $A$. Furthermore, if
$J\!\subseteq\!\{i_1,\ldots,i_{\ell+2}\}$ and
$\sum_{j\in J} a_j\!=\!\bO$, then $J\!=\!\{i_1,\ldots,i_{\ell+2}\}$. \qed
\end{prop}
We call $\Sigma$ the {\em unique non-degenerate sub-circuit} of $A$. 
Note that any $A\!=\!\{a_1,\ldots,a_{n+2}\}\!\subset\!\Zn$ with 
cardinality $n+2$, and $A$ not lying in any affine hyperplane, is a circuit. 
\begin{ex} 
It is easily checked that $A\!=$\scalebox{.7}[.7]{$\left\{
\begin{bmatrix} 0 \\ 0 \\ 0 \\ 0 \end{bmatrix},
\begin{bmatrix} 1 \\ 0 \\ 0 \\ 0 \end{bmatrix},
\begin{bmatrix} 2 \\ 0 \\ 0 \\ 0 \end{bmatrix},
\begin{bmatrix} 0 \\ 1 \\ 0 \\ 0 \end{bmatrix},
\begin{bmatrix} 0 \\ 0 \\ 1 \\ 0 \end{bmatrix},
\begin{bmatrix} 0 \\ 0 \\ 0 \\ 1 \end{bmatrix}\right\}$}$\subset\!\R^4$ is a
degenerate circuit, and that letting $\Sigma$ consist of the
first $3$ points of $A$ yields the unique non-degenerate sub-circuit
of $A$. In particular, $\Sigma$ has the same minimal
circuit relation (up to sign) as the non-degenerate circuit
$\{0,1,2\}$ in $\R^1$. \dia
\end{ex}
\begin{lemma} 
\label{lemma:gale} 
\scalebox{.96}[1]{Suppose $F\!=\!(f_1,\ldots,f_n)\!\in\!\Z\!\left[x^{\pm 1}_1,
\ldots,x^{\pm 1}_n\right]^n$ is an $(n+2)$-nomial $n\times n$ 
system}\linebreak 
supported on a circuit $A\!=\!\{a_1,\ldots,a_{n+2}\}$ with 
unique non-degenerate sub-circuit $\Sigma=$\linebreak   
$\{a_1,\ldots,a_m,a_{n+1},a_{n+2}\}$ and all coordinates of the $a_j$ have  
absolute value at most $d$. Suppose also that $F$ has coefficient matrix 
$[c_{i,j}]\!\in\!\{-H,\ldots,H\}^{n\times (n+2)}$ with all $n\times n$ 
sub-matrices of 
\scalebox{.7}[.7]{$\begin{bmatrix}c_{1,1} & \cdots & c_{1,n} & c_{1,n+2}\\ 
\vdots & \ddots & \vdots \\ 
c_{n,1} & \cdots & c_{n,n} & c_{n,n+2}\end{bmatrix}$} non-singular and  
\scalebox{.7}[.7]{$\begin{bmatrix}c_{i,n+1} & c_{i,n+2}\\ 
c_{j,n+1} & c_{j,n+2}\end{bmatrix}$} non-singular for all $i\!\neq\!j$. 
Then in time $O\!\left(n^{1+\omega}\log^2(ndH)\right)$, 
we can give either a true declaration that 
$F$ has no positive roots, or find rational numbers $\gamma_{1,0},
\gamma_{1,1}, \ldots, \gamma_{n+1,0}, \gamma_{n+1,1}$ and nonzero integers 
$b_1,\ldots,b_{m+1}$ (with $1\!\leq\!m\!\leq\!n$) satisfying the following 
properties:\\ 
1.\ The number of roots of the function 
$L(u)\!:=\!\sum^{m+1}_{i=1} b_i\log|\gamma_{i,1}u
+\gamma_{i,0}|$ in the open interval\linebreak  
\mbox{}\hspace{.5cm}$I\!:=\!\{u\!\in\!\R \; | \; 
\gamma_{i,1}u+\gamma_{i,0}\!>\!0 \text{ for all } 
i\!\in\!\{1,\ldots,n+1\}\}$ is exactly the number of positive\linebreak 
\mbox{}\hspace{.5cm}roots of $F$.\\  
2.\ $I$ is non-empty and, for each $i\!\in\!\{1,\ldots,n+1\}$, we have 
$\max\{\gamma_{i,1},\gamma_{i,0}\}\!>\!0$. \\   
3.\ $L$ is a non-constant real analytic function on $I$.\\  
4.\ We have height bounds $h(b_i)\!=\!O(n\log(nd))$ and 
$h(\gamma_{i,j})\!=\!O(n\log(nH))$ for all $i$. 
\end{lemma} 
\begin{ex} 
\label{ex:easy} 
Returning to Example \ref{ex:rigged}, one can easily apply Gauss-Jordan 
elimination to the underlying linear combinations of monomials, and then 
divide every equation by the\linebreak 
\scalebox{.92}[1]{last monomial $x^{166}_3 x^{68}_4 x^{343}_5$, to reduce 
$F\!=\!\bO$ to the following system having the same roots in $(\Rs)^5$:}\\  
\scalebox{.78}[.8]{\vbox{
\begin{eqnarray*} 
x^{36}_1 x^{194}_2 x^{-116}_3 x^{14}_4 x^{-283}_5 \hspace{10cm} 
  & = & 16384cx^{58}_1 x^{194}_2 x^{-142}_3 x^{-32}_4 x^{-318}_5 + \frac{1}{4}\\
\mbox{}\hspace{2.5cm}x^{76}_1 x^{240}_2 x^{-166}_3 x^{-27}_4 x^{-342}_5 
 \hspace{7.5cm} 
 & = & \hspace{.2cm}4096cx^{58}_1 x^{194}_2 x^{-142}_3 x^{-32}_4 x^{-318}_5 + 1 \\ 
\mbox{}\hspace{5cm}x^{74}_1 x^{179}_2 x^{-141}_3 x^{-68}_4 x^{-286}_5 
 \hspace{5cm} 
  & = & \hspace{.4cm}256cx^{58}_1 x^{194}_2 x^{-142}_3 x^{-32}_4 x^{-318}_5 
   + 1 \\ 
\mbox{}\hspace{7.5cm}x^{25}_1 x^{203}_2 x^{-122}_3 x^{-67}_4 x^{-343}_5 
 \hspace{2.5cm} 
  & = & \hspace{.6cm}16cx^{58}_1 x^{194}_2 x^{-142}_3 x^{-32}_4 x^{-318}_5 
   + 1 \\ 
\mbox{}\hspace{10cm}x^{20}_1 x^{167}_2 x^{-102}_3 x^{-56}_4 x^{-275}_5 
 \hspace{0cm} 
  & = & \hspace{1cm} cx^{58}_1 x^{194}_2 x^{-142}_3 x^{-32}_4 x^{-318}_5 + 1 
\end{eqnarray*}}}\\  
Note in particular that this new system clearly reveals that all the 
roots of $F$ in $(\Rs)^5$ (for our earlier chosen values of $c$) must in fact 
lie in $\R^5_+$, since the right-hand sides are all positive on $(\Rs)^5$.  
The underlying circuit relation for the exponent vectors above is the 
same as the circuit relation for the exponent vectors of $F$: 
$b\!=\!(-2,2,-2,2,-2,1,1)^\top$. Part of the proof of Lemma \ref{lemma:gale}, 
applied to our example here, will imply that the resulting linear combination 
of logarithms $L(u)$ can be easily read from $b$ and the righthand sides of 
our reduced system:\\ 
\scalebox{.9}[1]{$-2\log\left|16384cu+\frac{1}{4}\right|
+2\log\left|4096cu+1\right| -2\log\left|256cu+1\right| +2\log\left|16cu+1
\right| -2\log|cu+1| +\log|u|$.}\linebreak 
In particular, for any $c\!>\!0$, the number of roots of $L$ in 
$I\!=\!\R_+$ is the same as the number of roots of $F$ in $\R^5_+$. 
Our family of examples here is in fact an obfuscated version of a family 
derived in \cite{pr}, thus accounting for the nice coefficients and high 
number of positive roots ($6$) for $c\!\in\!\left[\frac{1}{20730},
\frac{1}{14392}\right]$. A more realistic example of coefficient growth can 
be found in Example \ref{ex:random} below (see also Example \ref{ex:random2} 
from Section \ref{sec:crit}). \dia 
\end{ex} 
\begin{ex} 
\label{ex:gen} 
The stated indexing of the non-degenerate sub-circuit $\Sigma$ within 
the union of supports $A$ is important, for otherwise it is 
possible for $F$ to have infinitely many positive roots. 
For instance, the $5$-nomial $3\times 3$ system  
$(x_2-1,x_1x_2 -x_3 -1, x^2_1x_2  -x_3 -1)$, with union of 
supports ordered so that $a_4\!=\!(0,0,1)^\top$ and $a_5\!=\!(0,0,0)^\top$,  
has $[c_{i,j}]^3_{i,j=1}$ non-singular (among other $3\times 3$ sub-matrices), 
but infinitely many positive roots: $(1+t,1,t)$ for all $t\!\in\!\R$. Note 
that, unlike the statement of Lemma \ref{lemma:gale}, the indexing is such 
that the last $2$ exponent vectors $\{a_4,a_5\}$ 
are {\em not} in the unique non-degenerate 
sub-circuit $\{(0,1,0)^\top,(1,1,0)^\top,(2,1,0)^\top\}$. \dia  
\end{ex} 

\noindent 
{\bf Proof of Lemma \ref{lemma:gale}:} We can divide 
$f_1,\ldots,f_n$ by $x^{a_{n+2}}$ without affecting the positive roots of $F$. 
So we may assume $a_{n+2}\!=\!\bO$ and, since this at worst doubles our 
original $d$, our $O$-estimates will be unaffected. We can then apply Lemma 
\ref{lemma:lin}, thanks to our assumption on the $n\times n$ sub-matrices  
of $[c_{i,j}]$, to reduce $F\!=\!\bO$ to a system of equations of the form 
$G\!=\!\bO$, {\em having the same solutions in $\Rn$ as} $F$, where 
$G\!:=\!(g_1,\ldots,g_n)$, 
\begin{eqnarray*} 
g_i(x)& =& x^{a_i}-\gamma_{i,1}x^{a_{n+1}}-\gamma_{i,0} \text{ for all } i,
\end{eqnarray*}  
\scalebox{.91}[1]{the $\gamma_{i,j}$ have logarithmic height 
$O(n\log(nH))$, and this reduction takes time just 
$O\!\left(n^{1+\omega} \log^2(nH)\right)$.}\linebreak  
To complete our notation, let us set $\gamma_{n+1,1}\!:=\!1$ and 
$\gamma_{n+1,0}\!:=\!0$.  

Clearly, if there is an $i$ such that both $\gamma_{i,0}$ and $\gamma_{i,1}$ 
are non-positive, then $G$ (and thus $F$) has no positive roots, and 
we can simply stop, having spent time just $O\!\left(n^{1+\omega}
\log^2(nH)\right)$. So we may assume the following: 
\begin{eqnarray}
\label{assu:pos} 
\text{For each } i\!\in\!\{1,\ldots,n+1\} \text{ we have } 
\max\{\gamma_{i,1},\gamma_{i,0}\}\!>\!0. 
\end{eqnarray} 
We can easily check whether $I$ is 
non-empty after sorting the (possibly infinite) numbers 
$-\gamma_{i,0}/\gamma_{i,1}$, using just $O(n\log n)$ comparisons of integers 
with $O(n\log(nH))$ bits (via, say, merge sort \cite{clrs}). 
If $I$ is empty then we can conclude that $F$ has no positive 
roots and stop (having spent time just $O\!\left(n^{1+\omega}\log^2(nH)
\right)$). So we may also assume the following: 
\begin{eqnarray}
\label{assu:nonempty} 
I \text{ is non-empty.} 
\end{eqnarray}  

We now establish Assertions (1)--(4) via $G$ and $\Sigma$:  
First let $b\!\in\!\Z^{(m+2)\times 1}$ be the unique (up to sign) 
minimal circuit relation of $\Sigma$. Note that the coordinates of $b$ are of 
logarithmic height $O(m\log(mH))$, and the computation of $b$ takes time 
$O\!\left(mn^{\omega}\log^2(nd)\right)$, thanks to Lemmata \ref{lemma:lin} and 
\ref{lemma:rightnull}.  

Observe that any root $\zeta\!\in\!\Rsn$ of $G$ must satisfy  
\begin{eqnarray} 
\label{eqn:ckt} 
\mbox{}\hspace{1cm}
\scalebox{.9}[1]{\text{$1 =  (\zeta^{a_1})^{b_1}\cdots
(\zeta^{a_m})^{b_m}(\zeta^{a_{n+1}})^{b_{m+1}} 
 =  (\gamma_{1,1}\zeta^{a_{n+1}}+\gamma_{1,0})^{b_1}
\cdots (\gamma_{m,1}\zeta^{a_{n+1}}+\gamma_{m,0})^{b_m}
(\zeta^{a_{n+1}})^{b_{m+1}} 1^{b_{m+2}}$,}}   
\end{eqnarray} 
So let $P(u)\!:=\!(\gamma_{1,1}u+\gamma_{1,0})^{b_1}
\cdots (\gamma_{m,1}u+\gamma_{m,0})^{b_m}u^{b_{m+1}}-1$. 
Note that $n\!=\!1\Longrightarrow (\gamma_{1,1},
\gamma_{1,0})\!=\!(-c_2/c_1,-c_3/c_1)\!\in\!(\Rs)^2$ and thus 
$P$ is a non-constant real rational function when $n\!=\!1$. So $L(u)$ is a 
non-constant real analytic function on $I$ when $n\!=\!1$. 
Let us then assume $n\!\geq\!2$. By Cramer's Rule, and our assumption on the 
$n\times n$ sub-matrices of $[c_{i,j}]$, we have $\gamma_{i,0}\!\neq\!0$ for 
all $i$. 
If $m\!=\!1$ then we have that $P$ is a non-constant real rational function 
since $b_{m+1}\!\neq\!0$ (thanks to $\Sigma$ being a non-degenerate circuit and 
Proposition \ref{prop:ckt}) and there is thus no way to cancel the 
$u^{b_{m+1}}$ factor in the product term of $P$. So $L(u)$ is a non-constant 
real analytic function on $I$ when $m\!=\!1$, and we thus now assume 
$m\!\geq\!2$. 

By our assumption on the rightmost 
$2\times 2$ minors of $[c_{i,j}]$, we see that $\gamma_{i,1}\!\neq\!0$ 
for some $i\!\in\!\{1,\ldots,m\}$, and thus $P$ must have 
$-\gamma_{i,0}/\gamma_{i,1}$ ($\neq\!0$) 
as a zero of possibly non-positive order. In particular, the order must be 
nonzero, thanks to Proposition \ref{prop:ckt}. 
So we see that $P$ is a\linebreak 
\scalebox{.97}[1]{non-constant real rational function and, again, 
$L$ is a non-constant real analytic function on $I$.}     

Now observe that any root $\zeta\!\in\!\Rn_+$ of $F$ yields  
$\zeta^{a_{n+1}}\!\in\!I$ as a root of $P$ by Equation (\ref{eqn:ckt}).
Moreover, by Proposition \ref{prop:mono}, {\em any} root 
$\zeta'\!\in\!\Rn_+$ of $F$ with $(\zeta')^{a_{n+1}}\!=\!\zeta^{a_{n+1}}$ 
must satisfy $\zeta'\!=\!\zeta$, since $G$ reduces to a binomial system 
with a unique positive root once the value of $x^{a_{n+1}}$ is fixed. 
(This is because the vectors $a_1,\ldots,a_n$ are linearly independent, 
thanks to $\{a_{n+1},\bO\}\!\subset\!\Sigma$ and $A$ being a circuit.) So $P$ 
has at least as many roots in $I$ as $F$ has in $\Rn_+$. 

Conversely, Proposition \ref{prop:mono} tells us that any root $u\!\in\!I$  
of $P$ yields a unique $\zeta\!\in\!\Rn_+$ satisfying $\left(\zeta^{a_1},
\ldots,\zeta^{a_n}\right)\!=\!(\gamma_{1,1}u+\gamma_{1,0},\ldots,\gamma_{n,1}u
+\gamma_{n,0})$. Recall that $b_{m+1}\!\neq\!0$. So we also obtain 
$\zeta^{a_{n+1}}=\left((\gamma_{1,1}u+\gamma_{1,0})^{-b_1}
\cdots(\gamma_{m,1}u+\gamma_{m,0})^{-b_m}\right)^{1/b_{m+1}}
\!=\!u^{b_{m+1}/b_{m+1}}\!=\!u$ by the definition of $P$. So $\zeta$ is 
in fact a root of $G$. Similarly, a root $u'\!\in\!I$ of $P$ with 
$u'\!\neq\!u$ would yield a positive root of $\left((\zeta')^{a_1},
\ldots,(\zeta')^{a_n}\right)\!=\!(\gamma_{1,1}u'+\gamma_{1,0},
\ldots,\gamma_{n,1}u' +\gamma_{n,0})$ with 
$(\zeta')^{a_{n+1}}\!\neq\!\zeta^{a_{n+1}}$ and thus a root 
$\zeta'\!\neq\!\zeta$ of $F$. So $F$ has at least as many roots in $\Rn_+$ as 
$P$ has in $I$. 

Observing that $P(u)\!=\!0 \Longleftrightarrow L(u)\!=\!0$ (for $u\!\in\!I$), 
and recalling Assumptions (\ref{assu:pos}) and (\ref{assu:nonempty}), 
we thus obtain Assertions (1)--(4). Noting that $m\!\leq\!n$, 
we are done. \qed 

\medskip 
Our sub-matrix condition from Lemma \ref{lemma:gale} in fact holds for a large 
fraction of integer coefficients: 
\begin{cor} 
\label{cor:gen} 
The fraction of matrices $[c_{i,j}]\!\in\!\{-H,\ldots,H\}^{n\times
(n+2)}$ with all $n\times n$ sub-matrices of
\scalebox{.7}[.7]{$\begin{bmatrix}c_{1,1} & \cdots & c_{1,n} & c_{1,n+2}\\ 
\vdots & \ddots & \vdots \\ 
c_{n,1} & \cdots & c_{n,n} & c_{n,n+2}\end{bmatrix}$} non-singular and 
\scalebox{.7}[.7]{$\begin{bmatrix}c_{i,n+1} & c_{i,n+2}\\ 
c_{j,n+1} & c_{j,n+2}\end{bmatrix}$} non-singular for all $i\!\neq\!j$ 
is at least $1-\frac{2n^2}{2H+1}$. Also, the fraction of 
matrices $[c_{i,j}]\!\in\!\{-H,\ldots,H\}^{n\times (n+1)}$ with 
leftmost $n\times n$ sub-matrix non-singular is at least $1-\frac{n}{2H+1}$. 
\end{cor} 

\noindent 
{\bf Proof:} Schwartz's Lemma \cite{schwartz} is a classic result 
that tells us that if $f\!\in\!\C[z_1,\ldots,z_n]$ has degree $d$ 
and $S\!\subset\!\C$ is a set of finite cardinality $N$, then 
$f$ vanishes at no more than $dN^{n-1}$ points of $S^n$. The 
condition stated in our corollary is simply the non-vanishing 
of a product of $n+1$ many $n\times n$ sub-determinants of 
$[c_{i,j}]$ and $\binom{n}{2}$ many $2\times 2$ sub-determinants 
of $[c_{i,j}]$. The resulting polynomial clearly has degree 
$(n+1)n+\frac{n(n-1)}{2}\cdot 2\!=\!2n^2$. Taking $S\!=\!\{-H,\ldots,H\}$ 
and applying Schwartz's Lemma, we obtain our first 
bound. Our second bound follows almost identically, just 
considering one determinant instead. \qed 

\medskip
Recall that a {\em critical point} of a function $L : \R\longrightarrow \R$ is
a root of the derivative $L'$.
\begin{prop}
\label{prop:crit}
Following the notation and assumptions of Lemma \ref{lemma:gale},
let $u_0\!:=\!\inf I$, $u_k\!:=\!\sup I$, and suppose
$u_1\!<\cdots<\!u_{k-1}$ are the critical points of $L$ in $I$
($k\!=\!1$ implying {\em no} critical points). Then the number of positive
roots of $F$ is exactly the number of\linebreak
$i\!\in\!\{0,\ldots,k-1\}$ such that
$\left(\lim_{u\rightarrow u^+_i} L(u)\right) 
\left(\lim_{u\rightarrow u^-_{i+1}} L(u)\right)\!<\!0$,
{\em plus} the number of degenerate roots of $L$ in $I$.
\end{prop}

\noindent
{\bf Proof:} It is clear that $L$ is strictly monotonic on any open
sub-interval $(u_i,u_{i+1})$ of $I$. So the image of $(u_i,u_{i+1})$
under $L$ is \\
\mbox{}\hfill $L_i\!:=\!\left(\min\left\{\lim_{u\rightarrow u^+_i} L(u),
\lim_{u\rightarrow u^-_{i+1}} L(u)\right\},
\max\left\{\lim_{u\rightarrow u^+_i} L(u),
\lim_{u\rightarrow u^-_{i+1}} L(u)\right\}\right)$, \hfill\mbox{}\\
and we see by the Intermediate Value Theorem that $L_i$ does {\em not}
contain $0 \Longleftrightarrow 
[\lim_{u\rightarrow u^+_i} L(u)$ and $\lim_{u\rightarrow u^-_{i+1}} L(u)$
are both non-positive or both non-negative$]$. So by Lemma \ref{lemma:gale},
we are done. \qed

\medskip 
We can now state analogues of Lemma \ref{lemma:gale} and Proposition \ref{prop:crit} 
for roots in $\Rsn$. 
\begin{lemma}
\label{lemma:gale+} 
Following the notation and assumptions of Lemma \ref{lemma:gale}, 
assume further that $b_{m+1}$ is odd, $a_{n+2}\!=\!\bO$, and let 
$\cA\!:=\![a_1,\ldots,a_n]$. Let $U\cA V\!=\!S$ be a Smith factorization 
of $\cA$, and $r$ the rank, over the field $\F_2$, of the mod $2$ reduction of 
$\cA$. Also, for any $u\!\in\!\R$, let $\eps_i\!:=\!\sign(\gamma_{i,1}u
+\gamma_{i,0})$, $\Lambda(u)\!:=\!\prod\limits^m_{i=1}\eps^{b_i\mathrm{mod} 
\; 2}_i$, and $(\Gamma'_1(u),\ldots,\Gamma'_n(u))\!:=\!(\eps_1,\ldots,
\eps_n)^{V \mathrm{mod} \; 2}$. Then the number of roots of $F$ in $\Rsn$ is 
exactly $2^{n-r}$ times the number of roots $u\!\in\!\R$ of $L$ satisfying  
both $\Lambda(u)\!=\!\sign(u)$ and $\Gamma'_{r+1}(u),\ldots,\Gamma'_n(u)\!>\!0$.
\end{lemma}

\begin{ex}
\label{ex:random} 
Consider the $6$-nomial $4\times 4$ system $F\!=\!(f_1,\ldots,f_4)$ defined by 
\begin{eqnarray*} 
\left( -12x^{8}_1  x^{18}_2          x^{16}_4 
       -5x^{4}_1  x_2      x^3_3    x^8_4 
      +17x^{11}_1 x^{19}_2 x_3      x^{17}_4 
       -4x^{11}_1 x^9_2    x^{14}_3 
       +2         x^{18}_2 x^{13}_3 x^{17}_4  
       +3x^5_1             x^{14}_3 x^{16}_4, \right. \\ 
-9x^{8}_1  x^{18}_2          x^{16}_4 
+14x^{4}_1  x_2      x^3_3    x^8_4  
-8x^{11}_1 x^{19}_2 x_3      x^{17}_4 
+3x^{11}_1 x^9_2    x^{14}_3 
+12         x^{18}_2 x^{13}_3 x^{17}_4 
- x^5_1             x^{14}_3 x^{16}_4, \\  
5x^{8}_1  x^{18}_2          x^{16}_4 
 +4x^{4}_1  x_2      x^3_3    x^8_4
+11x^{11}_1 x^{19}_2 x_3      x^{17}_4 
-16x^{11}_1 x^9_2    x^{14}_3 
+18         x^{18}_2 x^{13}_3 x^{17}_4 
-19x^5_1             x^{14}_3 x^{16}_4, \\  
\left. -x^{8}_1  x^{18}_2          x^{16}_4  
 +2x^{4}_1  x_2      x^3_3    x^8_4 
+11x^{11}_1 x^{19}_2 x_3      x^{17}_4 
-17x^{11}_1 x^9_2    x^{14}_3 
-14         x^{18}_2 x^{13}_3 x^{17}_4 
 -6x^5_1             x^{14}_3 x^{16}_4 \right). 
\end{eqnarray*} 
Proceeding as in Lemmata \ref{lemma:gale} and \ref{lemma:gale+}, we 
see that Gauss-Jordan elimination on the coefficient matrix, and 
computing the circuit relation underlying the exponent vectors, yields 
the following linear combination of logarithms:  
\begin{eqnarray*} 
L(u)&=&54667\log\left|\frac{39898}{27281}-\frac{84556}{27281}u\right|
      -16978 \log \left|\frac{47210}{27281}-\frac{125680}{27281}u\right| 
 \\ & & -43727 \log \left|\frac{42139}{27281}-\frac{126754}{27281}u\right| 
   +5123\log\left|\frac{20845}{27281}-\frac{114296}{27281}u\right|-10129\log|u|
\end{eqnarray*}
In particular, $L$ is an analytic function on\\ 
\mbox{}\hfill  
$\R\setminus\{0,0.182377...,0.332447...,0.375636...,0.471852...\}$\hfill
\mbox{}\\ 
whose roots encode the roots of $F$ in $(\Rs)^4$: Observing that\\ 
\mbox{}\hfill  
$\Lambda(u)\!=\!\sign\!\left(39898-84556u\right)\sign\!\left(42139
-126754u\right) \sign\!\left(20845-114296u\right)$,\hfill\mbox{}\\ 
we see that the only open 
intervals containing $u$ satisfying $\Lambda(u)\!=\!\sign(u)$ are\\ 
\mbox{}\hfill $(0,0.182377...)$, 
$(0.332447...,0.375636...)$, $(0.375636...,0.471852...)$. \hfill\mbox{}\\ 
(The mod $2$ reduction of our $\cA$ here has full rank $r\!=\!4$ and thus the 
condition involving the $\Gamma'_i(u)$ becomes vacuously true.) 
It is then easily checked that $L$ is strictly decreasing,  
with range $\R$, on the first and third intervals; and $L$ is positive on the 
second interval. (See also Corollary \ref{cor:crit} below.) So $L$ has exactly 
$2$ roots in $\Rs$ satisfying the sign conditions\footnote{$L$ also happens to 
be increasing, with range $\R$, on $(-\infty,0)$ and $(0.182377...,
0.332447...)$, and thus $L$ has $2$ more roots in $\Rs$ that do {\em not}  
satisfy the necessary sign conditions.} of Lemmata \ref{lemma:gale+}, and thus 
$F$ has exactly $2$ roots in $(\Rs)^4$. {\tt PHCpack} (a state of the art 
polyhedral homotopy solver \cite{verschelde}) confirms our root count for this 
$F$ in about 15 minutes, along with a count of $70834$ for the total number of 
roots in $(\Cs)^4$, as well as approximations of all these roots to $14$ 
decimal places. Our {\tt Maple} simulation counts the roots of $F$ in 
$(\Rs)^4$ and $\R^4_+$ in under one second. \dia  
\end{ex} 

\noindent 
{\bf Proof of Lemma \ref{lemma:gale+}:} Continuing the notation of the proof 
of Lemma \ref{lemma:gale}, we need to revisit the properties of the rational 
function $P$ defined earlier. In particular, whereas before we had a natural 
bijection between the roots of $F$ in $\Rn_+$ and the roots of $P$ in a 
particular interval $I$, we now need to consider roots of $F$ with negative 
coordinates and roots of $P$ throughout $\R$. In particular, a key difference 
from our earlier lemma is the following simple equivalence, valid for 
all $u\!\in\!\R$: $P(u)\!=\!0 \Longleftrightarrow [L(u)\!=\!0$ {\em and} 
$\Lambda(u)\!=\!  \sign(u)]$. (Indeed, we could obtain $u$ with $P(u)\!=\!-2$ 
without the condition involving $\Lambda(u)$.) Note also that by construction, 
$P(0)$ is either $1$ or undefined. 

So let $\zeta\!\in\!\Rsn$ be a root of $F$. 
By Relation (\ref{eqn:ckt}), $\zeta^{a_{n+1}}$ must be a nonzero real root of 
$P$ and, by the definition of $G$ and the $\gamma_{i,j}$ (and Proposition 
\ref{prop:rsn}), we must have $\Gamma'_{r+1}\!\left(\zeta^{a_{n+1}}\right),
\ldots,$\linebreak 
$\Gamma'_n\!\left(\zeta^{a_{n+1}}\right)\!>\!0$. By Proposition 
\ref{prop:rsn}, there must also be exactly $2^{n-r}$ many 
$\zeta'\!\in\!\Rsn$ of $F$ with $(\zeta')^{a_{n+1}}\!=\!\zeta^{a_{n+1}}$, 
because $G$ reduces to a binomial system once the value of $\zeta^{a_{n+1}}$ is 
fixed. So $F$ has no more than $2^{n-r}$ times as many roots in $\Rsn$ as $P$ 
has in $\Rs$. 

Conversely, if $u\!\in\!\Rs$ is a root of $P$ then Proposition \ref{prop:rsn} 
tells us that $\Gamma'_{r+1}(u), \ldots,\Gamma'_n(u)\!>\!0$ implies that there 
are exactly $2^{n-r}$ many $\zeta\!\in\!\Rsn$ satisfying\\ 
\mbox{}\hfill  
$(\zeta^{a_1},\ldots,\zeta^{a_n})\!=\!(\gamma_{1,1}u+\gamma_{1,0},\ldots,
\gamma_{n,1}u+\gamma_{n,0})$. \hfill\mbox{}\\ 
(Note also that $\gamma_{i,1}u+\gamma_{i,0}
\!\neq\!0$ for all $i$ since $P(u)\!\neq\!0$ when 
$\gamma_{i,1}u+\gamma_{i,0}\!=\!0$.) We 
then also obtain $\zeta^{b_{m+1}a_{n+1}}=(\gamma_{1,1}u+\gamma_{1,0})^{-b_1}
\cdots(\gamma_{m,1}u+\gamma_{m,0})^{-b_m} 
\!=\!u^{b_{m+1}}$ by Relation (\ref{eqn:ckt}). Since $b_{m+1}$ is odd, 
all our resulting $\zeta$ must satisfy $\zeta^{a_{n+1}}\!=\!u$ and 
therefore be roots of $G$ (and thus of $F$). Similarly, a real 
root $u'$ of $P$ with $(u')^{b_{m+1}}\!\neq\!u^{b_{m+1}}$ would 
yield a collection of $2^{n-r}$ many $\zeta'\!\in\!\Rsn$ that are roots of $F$ 
but with $(\zeta')^{a_{n+1}}\!\neq\!\zeta^{a_{n+1}}$, since 
$b_{m+1}$ is odd and $u\!\in\!\Rs$. So the number of roots of $F$ in $\Rsn$ 
is at least $2^{n-r}$ times the number of roots of $P$ in $\Rs$. \qed 

\medskip 
The following variant of Proposition \ref{prop:crit} can be proved 
almost the same as Proposition \ref{prop:crit}, simply using Lemma 
\ref{lemma:gale+} instead of Lemma \ref{lemma:gale}:  
\begin{cor} 
\label{cor:crit} 
Following the notation and assumptions of Lemma \ref{lemma:gale+},
let $w_1\!<\cdots<\!w_{\ell-1}$ be the critical points {\em and} 
poles of $L$ in $\R$, and set $w_0\!:=\!-\infty$ and 
$w_\ell\!:=\!+\infty$. Let $\cN$ be the number of 
$i\!\in\!\{0,\ldots,\ell-1\}$ such that 
$\Lambda(u)\!=\!\sign(u)$ and $\Gamma'_{r+1}(u),\ldots,\Gamma'_n(u)\!>\!0$ 
for all $u\!\in\!(w_i,w_{i+1})$, {\em and}   
$\left(\lim_{u\rightarrow w^+_i} L(u)\right) 
\left(\lim_{u\rightarrow w^-_{i+1}} L(u)\right)\!<\!0$. 
Then the number of roots of $F$ in $\Rsn$ is exactly $\cN$  
{\em plus} the number of degenerate roots of $L$ in $\R$. \qed 
\end{cor} 

\subsection{Heights of Algebraic Numbers and Linear Forms in Logarithms} 
\label{sub:linlogalg} 
Recall that if $\beta$ is in the algebraic closure $\overline{\Q}$ of $\Q$, 
with minimal polynomial $m(x_1):=c_0+\cdots+c_dx^d_1\!\in\!\Z[x_1]$ 
satisfying $\gcd(c_0,\ldots,c_d)\!=\!1$, then we may 
define the {\em (absolute) logarithmic height} of $\beta$ to be \\
\mbox{}\hfill $\displaystyle{h(\beta)\!:=\!\frac{1}{d}\left(
\log|c_d|+\sum\limits^d_{i=1}\log\max\{|\beta_i|,1\}\right)}$,\hfill\mbox{}\\ 
where $\beta_1,\ldots,\beta_d$ (among them, $\beta$) are all the roots of $m$. 
This definition in fact agrees with our earlier definition for rational 
numbers. Since $m$ must be irreducible we have 
$\#\{\beta_1,\ldots,\beta_d\}\!=\!d$. 
\begin{prop} 
\label{prop:heightsum} 
(See, e.g., \cite[Prop.\ 1.5.15, pg.\ 18]{bombgu}.) 
If $\alpha_1,\ldots,\alpha_k\!\in\!\overline{\Q}$ then   
$h\!\left(\sum\limits^k_{i=1}\alpha_i\right)$ is no greater than 
$\log(k)+\sum\limits^k_{i=1} h(\alpha_i)$. Also, 
$h\!\left(\prod\limits^k_{i=1}\alpha_i\right)\!\leq\!\sum\limits^k_{i=1} 
h(\alpha_i)$. \qed 
\end{prop} 

Letting 
$\left|c_0+c_1x_1+\cdots+c_dx^d_1\right|_2\!:=\!\sqrt{\sum^d_{i=0}|c_i|^2}$, 
we recall the following classical inequality: 
\begin{landau} 
\cite{mignotte} If $\beta\!\in\!\overline{\Q}$ has minimal polynomial 
$m\!\in\!\Z[x_1]$ with relatively prime coefficients then 
$\displaystyle{h(\beta)\!\leq\!\frac{\log|m|_2}{\deg m}}$. \qed 
\end{landau} 

It will also be useful to have a mildly refined version Liouville's classic 
bound \cite{liouville} 
on the separation between rational numbers and irrational algebraic numbers. 
\begin{thm} 
\label{thm:liouville} 
Suppose $\beta\!\in\!\overline{\Q}$, with
minimal polynomial $m\!\in\!\Z[x_1]$ of degree $d\!\geq\!2$. Then\\  
\mbox{}\hfill $\displaystyle{\left|\beta-\frac{p}{q}\right|<1\Longrightarrow 
\left|\beta-\frac{p}{q}\right|\geq 
\frac{\left(|m'(\beta)|+
\left|\frac{m''(\beta)}{2!}\right|+\cdots+\left|\frac{m^{(d)}(\beta)}{d!}
\right|\right)^{-1}}{q^d}}$\hfill\mbox{}\\ 
for all $p,q\!\in\!\Z$ with $q\!>\!0$.   
\end{thm} 

\noindent 
{\bf Proof:} First note that the numerator of the large fraction 
above is positive since $m$ is the minimal polynomial of $\beta$ 
and thus $m'(\beta)\!\neq\!0$.  

Via Taylor expansion we then obtain the following: 
\begin{eqnarray*}
|m(p/q)| & = & \left| m(\beta)+m'(\beta)\left(\frac{p}{q}-\beta\right)+ 
\frac{m''(\beta)}{2!}\left(\frac{p}{q}-\beta\right)^2+\cdots+
\frac{m^{(d)}(\beta)}{d!}\left(\frac{p}{q}-\beta\right)^d\right|\\ 
         & = & \left|\beta-\frac{p}{q}\right|\left|0
              +m'(\beta) + \frac{m''(\beta)}{2!}
              \left(\frac{p}{q}-\beta\right)+\cdots
          +\frac{m^{(d)}(\beta)}{d!}\left(\frac{p}{q}-\beta\right)^{d-1}\right|
\end{eqnarray*} 
\begin{eqnarray*} 
    &\leq & \left|\beta-\frac{p}{q}\right|\left(
            |m'(\beta)|+\left|\frac{m''(\beta)}{2!}\right|\left|\frac{p}{q} 
             -\beta\right|+\cdots +\left|\frac{m^{(d)}(\beta)}{d!}\right|
             \left|\frac{p}{q}-\beta\right|^{d-1}\right)\\ 
    & < & \left|\beta-\frac{p}{q}\right|\left(
          |m'(\beta)|+\left|\frac{m''(\beta)}{2!}\right|+\cdots
           +\left|\frac{m^{(d)}(\beta)}{d!}\right|\right). 
\end{eqnarray*} 

Since $m$ is irreducible and of degree $\geq\!2$, $m$ has no rational roots, 
and thus $q^d m(p/q)$ must be a nonzero integer. So we obtain 
$q^d|m(p/q)|\!\geq\!1$ and thus $|m(p/q)|\!\geq\!1/q^d$. Combined with 
our last Taylor series inequalities, we are done. \qed 

\medskip 
Finally, we recall the following paraphrase of a bound of  
Matveev \cite[Cor.\ 2.3]{matveev}, considerably strengthening earlier bounds of 
Baker and Wustholtz \cite{bakerwustholtz}. (See also \cite[Thm.\ 9.4]{bms}.) 
\begin{thm} 
\label{thm:bwm} 
Suppose $K$ is a degree $d$ real algebraic extension of $\Q$, 
$\alpha_1,\ldots,\alpha_m\!\in\!K\setminus\{0\}$, and 
$b_1,\ldots,b_m\!\in\!\Z\setminus\{0\}$. Let $B\!:=\!\max\{|b_1|,
\ldots,|b_m|\}$ and $\sA_i\!:=\!\max\{dh(\alpha_i),|\log \alpha_i|,0.16\}$ 
for all $i$. Then $\sum^m_{i=1} b_i\log\alpha_i\!\neq\!0$ implies that\\  
\mbox{}\hfill $\displaystyle{\log\left|
\sum^m_{i=1} b_i\log\alpha_i \right|>-
1.4\cdot m^{4.5} 30^{m+3}d^2(1+\log d)(1+\log(mB))\prod^m_{i=1}\sA_i}$.  
\hfill \qed 
\end{thm} 

\subsection{Bounds on Coefficients, Roots, and Derivatives of 
Univariate Polynomials} 
\label{sub:space} 
Letting $\left|c_0+c_1x_1+\cdots+c_dx^d_1\right|_\infty\!:=\!\max_i|c_i|$, 
recall the following classic bounds on the size and minimal spacing of roots of 
polynomials: 
\begin{prop} 
\label{prop:cauchy} 
(See, e.g., \cite[Thm.\ 8.1.4 \& Thm.\ 8.1.7, (i), (8.1.3)]{rs}.)  
If $f\!\in\!\Z[x_1]$ satisfies $|f|_\infty\!\leq\!H$ and 
$\zeta\!\in\!\C$ is a nonzero root of $f$ then 
$\frac{1}{1+H}\!<\!|\zeta|\!<\!1+H$. \qed 
\end{prop} 
\begin{mahler} \cite{mahler} 
Suppose $f\!\in\!\Z[x_1]$ is square-free, has degree $d$, and 
$|f|_\infty\!\leq\!H$. Then any two distinct complex roots $\zeta_1,\zeta_2$ 
of $f$ satisfy\\ 
\mbox{}\hfill $\displaystyle{|\zeta_1-\zeta_2|\!>\!\sqrt{3}(d+1)^{-(2d+1)/2}
H^{-(d-1)}}$. 
\hfill \mbox{}\\ 
In particular, $|\log|\zeta_1-\zeta_2||\!=\!O(d\log(dH))$. \qed 
\end{mahler} 
Letting $\left|c_0+c_1x_1+\cdots+c_dx^d_1\right|_1\!:=\!\sum^d_{i=0}|c_i|$, 
recall also the following nice bound on the coefficients of divisors of 
polynomials: 
\begin{lemma}
\label{lemma:factor}  
\cite[Thm.\ 4]{mignotte}
Suppose $f,g\!\in\!\C[x_1]$ have respective leading coefficients 
$c$ and $\gamma$, and $g|f$. Then 
$|g|_1\!\leq\!2^{\deg g}\left|\frac{\gamma}{c}\right||f|_2$. \qed 
\end{lemma} 

A simple consequence of the last two bounds is the following 
extension of Mahler's bound to the case of polynomials with degenerate roots. 
\begin{cor} 
\label{cor:rootspacing} 
Suppose $f\!\in\!\Z[x_1]$ has degree $d$ and $|f|_\infty\!\leq\!H$. 
Then any two distinct complex roots 
$\zeta_1,\zeta_2$ of $f$ satisfy $|\log|\zeta_1-\zeta_2||\!=\!
O(d^2+d\log(dH))$. \qed  
\end{cor} 

We will also need the following bound on the coefficients of products of 
polynomials: 
\begin{gelfond} (See, e.g., \cite[Lemma 1.6.11, pg.\ 27]{bombgu})  
(Special case) If $f_1,\ldots,f_k\!\in\Z[x_1]$ then 
$\left|\prod\limits^k_{i=1}f_i\right|_\infty\!\leq\!
\prod\limits^k_{i=1} 2^{\deg f_i}|f_i|_\infty$. \qed   
\end{gelfond} 

Finally, we will need the following bound on higher derivatives of 
polynomials, dating back to work of Duffin and Schaeffer \cite{ds}, based on a 
classic bound of A.\ A.\ Markov \cite{markov}: 
\begin{cor}
\label{cor:markov} 
Suppose $f\!\in\!\C[x_1]$ has degree $d$ and $t\!>\!0$. Then\\ 
\mbox{}\hfill $\displaystyle{\underset{-t\leq x_1 \leq t}{\max} 
\left|f^{(j)}(x_1)\right|\leq \frac{d^2(d^2-1^2)\cdots (d^2-(j-1)^2)}{1\cdot 3 
\cdots (2j-1)}\cdot \frac{\underset{-t\leq x_1 \leq t}{\max}|f(x_1)|}{t^j}}$. 
\hfill \qed 
\end{cor} 

\noindent 
After rescaling the variable so it ranges over $[-1,1]$, 
the statement above follows immediately from \cite[Thm.\ 15.2.6 \& 
Cor.\ 15.2.7, Sec. 15.2]{rs}. The latter results in fact 
include conditions for equality in the bound above. 

\section{Critical Values of Linear Forms in Logarithms and 
Their Signs} 
\label{sec:crit} 
We are now ready to prove two key lemmata (\ref{lemma:crit} and 
\ref{lemma:spacing} below) that enable our new complexity bounds. 
\begin{lemma}
\label{lemma:crit}
Suppose $m\!\geq\!2$, $b_i\!\in\!\Z\setminus\{0\}$ and
$\gamma_{i,1},\gamma_{i,0}\!\in\!\Q$ for all $i\!\in\!\{1,\ldots,m\}$, 
$\det$\scalebox{.7}[.7]{$\begin{bmatrix}\gamma_{i,1} & 
\gamma_{i,0}\\ \gamma_{j,1} & \gamma_{j,0}\end{bmatrix}$}$\neq\!0$\linebreak 
\scalebox{.95}[1]{for all $i$, $j$, $B\!:=\!\max_i|b_i|$, and $
h(\gamma_{i,j})\!\leq\!\log H$ for some $H\!>\!1$. 
Let $\displaystyle{L(u):=\sum^m_{i=1}b_i\log|\gamma_{i,1}u+\gamma_{i,0}|}$.}
\linebreak 
Then the critical points of $L$ in $\R$ are exactly the real roots of a 
polynomial $g\!\in\!\Z[u]$ of degree at most $m-1$ with   
$|g|_\infty\!\leq\!m2^{m-1}BH^{2m}$. In particular, 
$\log|g|_\infty\!=\!O(\log(B)+m\log H)$ and $L$ has at most $m$ real roots.  
\end{lemma} 

\begin{ex}
\label{ex:random2} 
Example \ref{ex:random} is more representative (than 
Example \ref{ex:easy}) of the coefficient growth one 
encounters when converting $F$ into a univariate linear combination of 
logarithms $L$: There we saw an input $6$-nomial $4\times 4$ system $F$ with
coefficients and exponents having at most $2$ digits, resulting in an
$L$ with coefficients having $6$ or fewer digits.
In particular, the polynomial $g$ encoding the critical points of $L$ is\\
\scalebox{.96}[1]{\vbox{
\begin{eqnarray*}
-85015812446550320118784u^4
+160578806134338659719072u^3
-78932164016242868100268u^2\\
+13833463598904597755876u
-837930167824219163155,
\end{eqnarray*}}}\\
which has coefficients with at most $24$ digits, and $2$ real roots, neither
of which lies in the sub-intervals of $\R$ contributing to the root count of
$F$ in $(\Rs)^4$. More to the point, it is the signs of the poles of 
$L$, instead of the signs of the critical values, that determine the 
number of roots of $F$ in $(\Rs)^4$ for Example \ref{ex:random}. \dia
\end{ex} 
\begin{ex}  
Returning to Example \ref{ex:easy}, which had $L(u)$ being\\
\scalebox{.9}[1]{$-2\log\left|16384cu+\frac{1}{4}\right|
+2\log\left|4096cu+1\right| -2\log\left|256cu+1\right| +2\log\left|16cu+1
\right| -2\log|cu+1| +\log|u|$,}\linebreak  
it is easily checked via {\tt Maple} that this $L$ has exactly $5$ critical 
values, alternating in sign, and the underlying critical points interlace the 
$6$ positive roots of $L$. \dia 
\end{ex}  

\noindent
{\bf Proof of Lemma \ref{lemma:crit}:} First observe that
$L'(u)\!=\!\sum\limits^m_{i=1} \frac{b_i\gamma_{i,1}}
{\gamma_{i,1}u +\gamma_{i,0}}$. Thanks to our determinant assumption,  
$L'$ has exactly $m$ distinct poles. Gelfond's Lemma implies 
that the coefficients of the monomial term expansion of 
$(r_1u+s_1)\cdots(r_ku+s_k)$ have logarithmic height $O(k\log H)$ if the 
$r_i$ and $s_i$ are integers of absolute value at most $H$. 
Letting $\nu_i$ be the least common multiple of the denominators of 
$\gamma_{i,1}$ and $\gamma_{i,0}$, and setting\\  
\mbox{}\hfill $\displaystyle{g_i(u)\!:=\!\left.\left(b_i\gamma_{i,1}\nu_i 
\prod\limits^m_{j=1}(\gamma_{j,1}u +\gamma_{j,0})\nu_j\right) 
\right/((\gamma_{i,1}u +\gamma_{i,0})\nu_i)}$,\hfill\mbox{}\\   
let us define $g(u)\!:=\!\sum^m_{i=1}g_i(u)$. Clearly, 
$g_i\!\in\!\Z[u]$ for all $i$, $g\!\in\!\Z[u]$, and 
$g(u)$ is nothing more than $L'(u)\prod^m_{j=1}
(\gamma_{j,1}u+\gamma_{j,0})\nu_j$. 
So we clearly obtain the statement on the real critical points of $L$ 
being the real roots of $g$, and it is clear that $\deg g\!\leq\!m-1$. 

\scalebox{.95}[1]{Gelfond's Lemma implies that $|g_i|_\infty\!\leq\!
BH^2\left(2^{m-1}H^{2(m-1)}\right)$. 
Clearly then, $|g|_\infty\!\leq\!m2^{m-1}BH^{2m}$.} 

\scalebox{.95}[1]{That $L$ has at most $m$ roots in $I$ is immediate from 
Rolle's Theorem, since $\deg g\!\leq\!m-1$. \qed}  

\medskip
\scalebox{.95}[1]{Recall that a {\em critical value} of a function 
$L : \R\longrightarrow \R$ is the value of $L$ at a critical point of $L$.}  
\begin{lemma}
\label{lemma:spacing} 
Following the notation and assumptions of Lemma \ref{lemma:crit}, 
let $I$ be any open interval defined by consecutive real poles of $L$,  
let $\eps$ denote any nonzero 
critical value of $L$, and let $\delta$ (resp.\ $\eta$) be the 
minimum of $|\zeta_1-\zeta_2|$ over all distinct roots 
$\zeta_1,\zeta_2\!\in\!I$ of $L$ (resp.\ the derivative $L'$). Finally, let 
$\Delta$ denote the minimum of $|\zeta-\mu|$ as $\zeta$ (resp.\ $\mu$) 
ranges over the critical points (resp.\ poles) of $L$. Then: \\ 
\mbox{}\hspace{.5cm}1.\ $\log \eta >-O(m\log(B)+m^2\log H)$,\\ 
\mbox{}\hspace{.5cm}2.\ $\log|\eps| > 
  -O\!\left(61^m\log^{m+1}\left(BH^{3m-1}\right)\right)$,\\ 
\mbox{}\hspace{.5cm}3.\ $\log\Delta > -O\!\left(m\log(B)+m^2\log H\right)$, 
and\\  
\mbox{}\hspace{.5cm}4.\ $\log \delta > 
 -O\!\left(61^m\log^{m+1}\left(BH^{3m-1}\right)\right)$. 
\end{lemma} 

\noindent
{\bf Proof:}  
If $L$ has no critical points then, by Rolle's Theorem, $L$ 
has at most $1$ root in $I$ and Assertions (1)--(3) are vacuously 
true. So let us assume $L$ has exactly $k-1$ critical points 
(with $k\!\geq\!2$) in the open
interval $I$, $u_0\!:=\!\inf I$, $u_k\!:=\!\sup I$, and suppose
$u_1\!<\cdots<\!u_{k-1}$ are the critical points of $L$ in $I$. 
Also let $g$ denote the polynomial from Lemma \ref{lemma:crit}. 
Below we illustrate a coarse approximation of what the graph of $L$ 
can look like, along with some of our notation:\\   
\mbox{}\hfill 
\scalebox{.8}[.65]{
\begin{picture}(470,190)(0,0)
\put(30,0){\epsfig{file=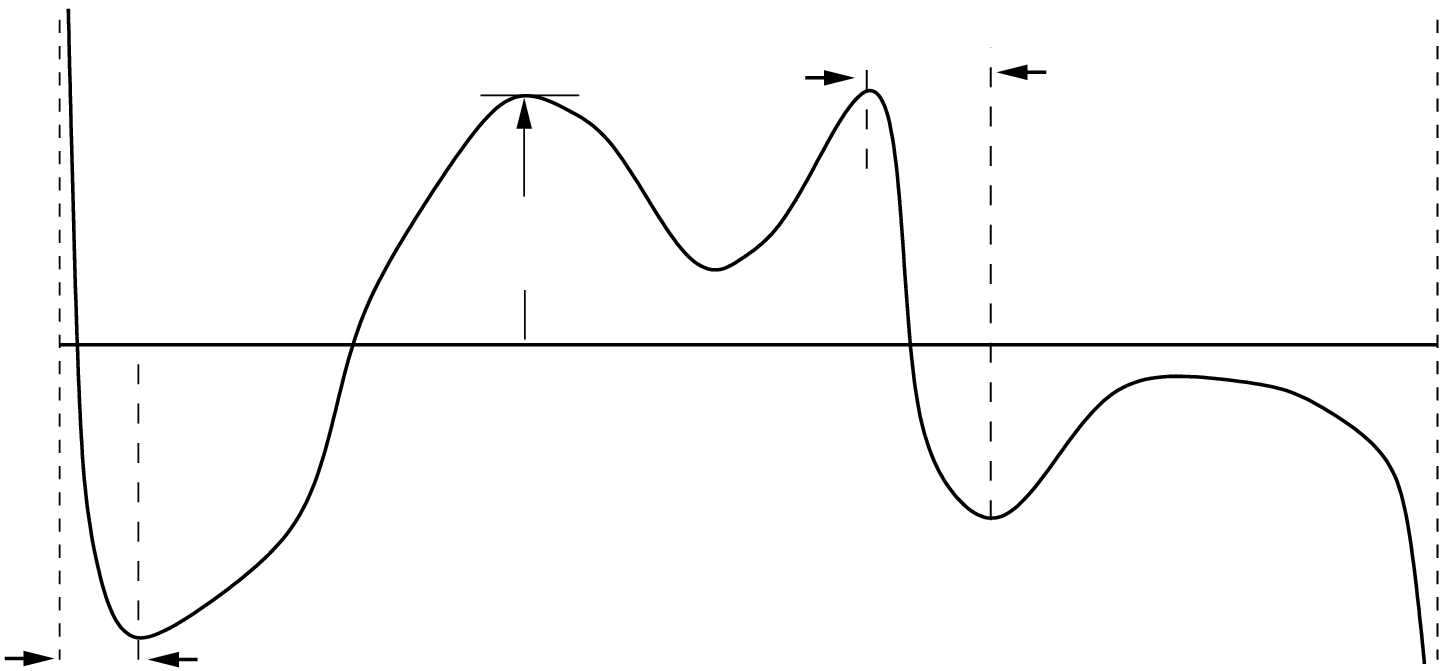,height=2.5in,clip=}} 
\put(46,88){\circle*{4}}
\put(33,82){\scalebox{1}[1]{$u_0$}}
\put(50,-3){\scalebox{1.4}[1.4]{$\Delta$}}
\put(62,91){\scalebox{1}[1]{$u_1$}}
\put(67.5,88){\circle*{4}}
\put(276,156){\scalebox{1.4}[1.4]{$\eta$}}
\put(168,112){\scalebox{1.4}[1.4]{$\eps$}}
\put(167,79){\scalebox{1}[1]{$u_2$}}
\put(171.5,88){\circle*{4}} 
\put(298,91){\scalebox{1}[1]{$u_5$}}
\put(297,88){\circle*{4}} 
\put(419,82){\scalebox{1}[1]{$u_7$}}
\put(417.5,88){\circle*{4}} 
\end{picture}}\hfill\mbox{} 

\scalebox{.95}[1]{Assertion (1) then follows immediately by applying Corollary 
\ref{cor:rootspacing} to $g$, thanks to Lemma \ref{lemma:crit}.}   

Assertion (2) then follows routinely from Theorem \ref{thm:bwm} upon observing 
that $|\eps|$ is nothing more than the absolute value of a linear 
combination of logarithms of real algebraic numbers. In particular, 
the arguments of the logarithms constituting $L(u_j)$ at a critical 
value $u_j\!\in\!I$ (for some $j\!\in\!\{1,\ldots,k-1\}$) all lie in the same 
real algebraic extension: $\Q(u_j)$. Noting that 
the minimal polynomial, $p$, of $u_j$ has degree $\leq\!m-1$, Lemmata 
\ref{lemma:factor} and \ref{lemma:crit} then tell us that 
$|p|_\infty\!\leq\!2^{m-1}|g|_\infty|g|_2$ (since $p|g$), and thus 
\begin{eqnarray}
\label{ineq:pinfinity} 
|p|_\infty \leq 2^{m-1}(mB2^{m-1}H^{2m})(\sqrt{m}(mB2^{m-1}H^{2m})) 
   =  m^{3/2}8^{m-1}B^2 H^{4m}. 
\end{eqnarray} 
So $\log|p|_2\!\leq\!\log\left(\sqrt{m}\cdot m^{3/2}B^2 8^{m-1}H^{4m}\right)$,  
and thus Landau's Inequality tells us that\\  
\mbox{}\hfill $\deg(p)h(u_j) \leq \log\left(m^28^{m-1} B^2 H^{4m}\right)$.
\hfill\mbox{}\\  
Proposition \ref{prop:cauchy} then tells us that 
$|u_j|< 1+|g|_\infty \leq 1+m2^{m-1}B H^{2m}$.  
Also, \linebreak 
$h(\gamma_{i,1}u_j+\gamma_{i,0})\!\leq\!\log(2)+(h(\gamma_{i,1})+
h(u_j))+h(\gamma_{i,0})$, thanks to Proposition \ref{prop:heightsum}. 
So then,\\  
\mbox{}\hspace{1cm}$\deg(p)h(\gamma_{i,1}u_j+\gamma_{i,0}) \leq 
(m-1)\log(2)+2(m-1)\log(H)+\log\!\left(m^2 8^{m-1}B^2H^{4m}\right)$ \\ 
\mbox{}\hspace{4.85cm}$\leq\log\!\left(m^2 16^{m-1}B^2H^{6m-2}\right)$.\\ 
Theorem \ref{thm:bwm} then tells us that 
\begin{eqnarray*}
\log |\eps| & > & \text{\scalebox{1}[1]{$-1.4\cdot 30^{m+3}m^{6.5}(1+\log(m-1))
 (1+\log(mB))\left(\log\!\left(m^2 16^{m-1}B^2 H^{6m-2}\right)\right)^m$}} \\ 
           & = & \text{\scalebox{1}[1]{$-1.4\cdot 30^3\cdot 30^m m^{6.5}
(1+\log(m-1))(1+\log(mB))\log^m\!\left(m^2 16^{m-1}B^2 H^{6m-2}\right)$}}\\  
           & = & -O\!\left( 30^m m^{6.5}(1+\log(m-1))(1+\log(mB)) 
\log^m\!\left(m^2 16^{m-1}B^2 H^{6m-2}\right) \right)\\ 
           & = & -O\!\left( 30^m m^{6.5}(2\log(m-1))(2\log(mB)) 
\log^m\!\left(m^2 16^{m-1}B^2 H^{6m-2}\right) \right)\\ 
           & = & -O\!\left( 30^m m^{6.5}\log(m-1)\log(mB) 
\log^m\!\left(m^2 16^{m-1}B^2 H^{6m-2}\right) \right).  
\end{eqnarray*} 
Since $m^{6.5}\log(m-1)\!=\!O\!\left(30^{\theta_1m}\right)$ 
and $m^2 16^{m-1}\!=\!O\!\left((B^2H^{6m-2})^{\theta_2}\right)$ for 
any $\theta_1,\theta_2\!>\!0$, we then obtain 
$\displaystyle{\log|\eps|  >  -O\!\left( 30^{(1+\theta_1)m} \log(mB) 
\log^m\!\left(\left(B^2 H^{6m-2}\right)^{1+\theta_2}\right) \right) }$    
\begin{eqnarray*} 
\mbox{}\hspace{1.2cm}  = & -O\!\left( 30^{(1+\theta_1)m} 
\log^{m+1}\!\left(\left(B^2H^{6m-2}\right)^{1+\theta_2}\right)\right) 
  \hspace{1.2cm} \mbox{} & \\  
           = & -O\!\left( 30^{(1+\theta_1)m}(1+\theta_2)^{m+1}2^{m+1} 
\log^{m+1}\!\left(BH^{3m-1}\right)\right) &   
           = -O\!\left( 61^m \log^{m+1}\!\left(BH^{3m-1}\right)\right),  
\end{eqnarray*} 
(picking, say, $\theta_1\!=\!\log(30.25/30)$ and $\theta_2\!=\!0.25/30.25$), 
thus proving Assertion (2). Note that we could have used any 
number strictly greater than $60$, instead of $61$, but we have opted for a 
more concise bound and a shorter proof. 

To prove Assertion (3) observe that 
$\Delta\!=\!\min\{u_1-u_0,u_k-u_{k-1}\}$, i.e.,  
$\Delta\!=\!|u_{j'}-u_\ell|$ for some $j'\!\in\!\{1,k-1\}$ and 
$\ell\!\in\!\{0,k-1\}$, by our earlier definitions. If 
$\Delta\!=\!\infty$ then there is nothing to prove, so let us 
assume $\Delta\!<\!\infty$. If $u_{j'}\!\in\!\Q$ then Assertion (3) follows 
easily from Proposition \ref{prop:heightsum}, since $u_\ell$ has logarithmic 
height no greater than $2\log H$ and $u_{j'}$ must have logarithmic height 
no greater than $O(\log(B)+m\log H)$. 
So we may assume that $u_{j'}$ is algebraic of degree at least $2$ over $\Q$.  

We can then apply Theorem \ref{thm:liouville} and Lemma \ref{lemma:crit} 
to obtain that $\Delta$ must be bounded from below by 
\begin{eqnarray}
\label{ineq:liouville} 
\frac{\left(|p'(u_{j'})|+
\left|\frac{p''(u_{j'})}{2!}\right|+\cdots+\left|\frac{p^{(m-1)}(u_{j'})}
{(m-1)!}\right|\right)^{-1}}{H^{2(m-1)}}.  
\end{eqnarray} 
We know that $|u_{j'}|< 1+m2^{m-1}B H^{2m}$ by Proposition \ref{prop:cauchy}, 
so it is enough to minimize the preceding sum of derivative norms over 
the interval $J\!:=\!\left[-1-m2^{m-1}B H^{2m},1+m2^{m-1}B H^{2m}\right]$.  

Noting the easy inequality 
$\displaystyle{\underset{-t\leq x_1 \leq t}{\max} 
\left|f(x_1)\right|\leq |f|_1\max\left\{1,|t|^d\right\}}$ for any 
$f\!\in\!\R[x_1]$ of degree $d$ and $t\!\in\!\R$ we then have  
\begin{eqnarray} 
\label{ineq:easy} 
\underset{x_1 \in J}{\max}|g(x_1)| & \leq & |g|_1\underset{x_1 \in J}
                                                   {\max}\max\left\{1,
                                                 x^{m-1}_1\right\} \nonumber \\ 
         & \leq & m|g|_\infty\left(1+m2^{m-1}BH^{2m}\right)^{m-1}, \nonumber \\ 
         & \leq & m\cdot m2^{m-1}BH^{2m}\left(1+m2^{m-1}BH^{2m}\right)^{m-1},
\end{eqnarray}  
and Corollary \ref{cor:markov} then implies  
\begin{eqnarray*}
|p^{(r)}(u_{j'})|& \leq & \underset{x_1 \in J }{\max} 
\left|p^{(r)}(x_1)\right| \\ 
& \leq &  \frac{(m-1)^2((m-1)^2-1^2)\cdots ((m-1)^2-(r-1)^2)}{1\cdot 3 
\cdots (2r-1)}\cdot |p|_1\cdot \left(1+m2^{m-1}BH^{2m}\right)^{m-1}, \\
\end{eqnarray*} 
since $p$ has degree $\leq\!m-1$. Since $|p|_1\!\leq\!m|p|_\infty$, 
Inequality (\ref{ineq:pinfinity}) tells us that  
$\left|p^{(r)}(u_{j'})\right|$ is bounded from above by\\  
\scalebox{.95}[1]{\vbox{\begin{eqnarray*} 
\frac{(m-1)^2((m-1)^2-1^2)\cdots ((m-1)^2-(r-1)^2)}{1\cdot 3 
\cdots (2r-1)}\cdot m \cdot \left(m^{3/2}8^{m-1}B^2 H^{4m}\right)
\left(1+m2^{m-1} B H^{2m}\right)^{m-1}  
\end{eqnarray*}}} 
\begin{eqnarray*} 
& \leq & \binom{(m-1)^2}{r} m \cdot \left(m^{3/2}8^{m-1} B^2 H^{4m}\right)
\left(1+m2^{m-1} B H^{2m}\right)^{m-1}  \\
& \leq & \left(\frac{(m-1)^2e}{r}\right)^r m^{5/2}8^{m-1}B^2 H^{4m} 
   \left(1+m2^{m-1} B H^{2m}\right)^{m-1}, \\
\end{eqnarray*}  
where the last inequality follows easily from Stirling's classical 
estimate for the factorial function. So then,  
$|p'(u_{j'})|+
\left|\frac{p''(u_{j'})}{2!}\right|+\cdots+\left|\frac{p^{(m-1)}(u_{j'})}
{(m-1)!}\right|$ is strictly less than \\ 
\mbox{}\hspace{.5cm}$\left[(m-1)^2e+\frac{((m-1)^2e/2)^2}{2!}+\cdots
+\frac{((m-1)^2e/r)^r}{r!}\right] m^{5/2}8^{m-1} B^2 H^{4m} 
\left(1+m2^{m-1} B H^{2m}\right)^{m-1}$. 

Now, by the Maclaurin series for 
$e^x$, the bracketed factor above is strictly less than $e^{(m-1)^2e}$. 
So then, $|p'(u_{j'})|+
\left|\frac{p''(u_{j'})}{2!}\right|+\cdots+\left|\frac{p^{(m-1)}(u_{j'})}
{(m-1)!}\right|$ is strictly less than  
\begin{eqnarray*} 
& & e^{(m-1)^2e} m^{5/2}8^{m-1} B^2 H^{4m} 
\left(1+m2^{m-1} B H^{2m}\right)^{m-1}\\ 
&=&O\!\left(e^{m(m-1)e} m^{5/2}8^{m-1} B^2 H^{4m}
\left(1+m2^{m-1} B H^{2m}\right)^{m-1} 
\right)\\  
&=&O\!\left(\left(e^{e+\theta_3}\right)^{m(m-1)} B^2 H^{4m} 
\left(2^{(1+\theta_4)(m-1)} B H^{2m}\right)^m\right)\\  
&=&O\!\left(\left(2^{1+\theta_4}e^{e+\theta_3}\right)^{m(m-1)}
B^{m+2}H^{2m^2+4m}\right), 
\end{eqnarray*} 
for any $\theta_3,\theta_4\!>\!0$. Observing that $2e^e\!<\!30.31$, 
we can then clearly pick $\theta_3$ and $\theta_4$ to obtain 
\begin{eqnarray*}
|p'(u_{j'})|+
\left|\frac{p''(u_{j'})}{2!}\right|+\cdots+\left|\frac{p^{(m-1)}(u_{j'})}
{(m-1)!} \right| & = & 
O\!\left(31^{m(m-1)}B^{m+2}H^{2m^2+4m}\right), 
\end{eqnarray*} 
and thus, combining with Inequality (\ref{ineq:liouville}), 
$\frac{1}{\Delta}\!=\!O\!\left(31^{m(m-1)}B^{m+2}H^{2m^2+6m-2}\right)$, 
and we obtain Assertion (3) by taking logarithms. 

To prove Assertion (4) we merely use the Mean Value Theorem. First, 
let $\delta'$ be the minimum distance between $u_j$ and $\zeta'$ 
where $\zeta'$ is any root of $L$. (Recall that $u_j$, with $j\!\in\!\{1,
\ldots,k-1\}$, is a critical point of $L$ corresponding to a nonzero 
critical value $\eps$ of $L$.)  
Clearly, $\delta\!\geq\!2\delta'\!>\!0$ (thanks to Rolle's Theorem), so 
it is enough to prove a 
sufficiently good lower bound on $\delta'$. Note in particular that if 
$\delta'\!>\!\Delta/2$ then we are done, thanks to Assertion (3). So 
let us assume $\delta'\!\leq\!\Delta/2$. 

Recall that, from the proof of Lemma \ref{lemma:crit}, we have  
$L'(u)\!=\!g(u)/\prod\limits^m_{i=1} (\gamma_{i,1}u+\gamma_{i,0})\nu_i$ where 
$\nu_i$ is the least common multiple of the denominators of
$\gamma_{i,1}$ and $\gamma_{i,0}$. 

Clearly, if $\gamma_{i,1}\!\neq\!0$, then  
$|(\gamma_{i,1}u+\gamma_{i,0})\nu_i|\!=\!
\left|u-\frac{-\gamma_{i,0}}{\gamma_{i,1}}\right||\gamma_{i,1}\nu_i|\!>\!
|\gamma_{i,1}\nu_i|\Delta/2\!\geq\!\Delta/2$ for all 
$u\!\in\!\left[u_j-\delta',u_j+\delta'\right]$, since $\zeta'\!\in\!I$, 
$\frac{-\gamma_{i,0}}{\gamma_{i,1}}$ is a pole of $L$, and 
$\gamma_{i,1}\nu_i$ is a nonzero integer. On the 
other hand, if $\gamma_{i,1}\!=\!0$, then 
$|(\gamma_{i,1}u+\gamma_{i,0})\nu_i|\!=\!|\gamma_{i,0}\nu_i|\!\geq\!1$ 
since $\gamma_{i,0}\nu_i$ is a nonzero integer. So then,  
\begin{eqnarray} 
\label{ineq:prod} 
\prod\limits^m_{i=1}(\gamma_{i,1}u+\gamma_{i,0})\nu_i>\Delta^m/2^m 
\text{ for all $u\!\in\!\left[u_j-\delta',u_j+\delta'\right]$.} 
\end{eqnarray}  

By the Mean Value Theorem, we must have 
$|L'(\xi)|\!=\left|\frac{\eps}{\delta'}\right|$ for some 
$\xi\!\in\!\left(u_j-\delta',u_j+\delta'\right)$. So then, thanks to 
Inequalities (\ref{ineq:easy}) and (\ref{ineq:prod}),  
we obtain  
\begin{eqnarray*}
|L'(\xi)|&=&\left|g(\xi)\left/
\prod\limits^m_{i=1}(\gamma_{i,1}\xi+\gamma_{i,0})\nu_i\right.\right|\\ 
& < & m^2 2^{m-1}BH^{2m}\left(1+m2^{m-1}BH^{2m}\right)^{m-1} 
\frac{2^m}{\Delta^m}\\ 
& \leq & m^2 2^{2m-1}BH^{2m}\left(1+m2^{m-1}BH^{2m}\right)^{m-1} 
O\!\left(31^{m(m-1)}B^{m+2} H^{2m^2+6m-2}\right)^m 
\end{eqnarray*} 
\scalebox{.95}[1]{Since $\delta'\!=\!|\eps/L'(\xi)|$ we thus obtain that 
$\log \delta' = \log|\eps|-\log|L'(\xi)|$, which is then 
bounded from}\linebreak 
\scalebox{.91}[1]{below by 
$-O\!\left(61^m\log^{m+1}\left(BH^{3m-1}\right)\right)
-O\!\left(\log(m)+m+\log(B)+m\log(H)
+m\log\!\left(m2^{m-1}BH^{2m} \right)\right)$}\linebreak 
\mbox{}\hspace{2cm}$-mO\!\left(m^2\log(31)+m\log(B)+m^2\log H\right)$,\\    
which reduces to 
$-O\!\left(61^m\log^{m+1}\left(BH^{3m-1}\right)\right)$. \qed

\subsection{The Complexity of Approximating Logarithms and Real Roots 
of Polynomials} 
\label{sub:cxity} 
Any real number can be expressed in binary. Since 
$2^{\floor{\log_2 x}}\!\leq\!x\!\leq\!2^{1+\floor{\log_2 x}}$ 
for any $x\!\in\!\R_+$, it is easy to check that $1+\floor{\log_2 x}$ is the 
number of bits for the integer part of $x$. It then makes sense to call the 
binary expansion of $\floor{2^{\ell-1-\floor{\log_2 x}}x}$ the 
{\em $\ell$ most significant bits} of an $x\!\in\!\R_+$. Clearly, knowing 
the $\ell$ most significant bits of $x$ means that one knows 
$x$ within a multiple of $(1+2^{-\ell})^{\pm 1}$. 

Let us recall the following classical fact on approximating 
logarithms via Arithmetic-Geometric Iteration: 
\begin{thm} 
\label{thm:log} 
\cite[Sec.\ 5]{dan}  
Given any positive $x\!\in\!\Q$ of logarithmic height $h$, and 
$\ell\!\in\!\N$ with $\ell\!\geq\!h$, we can compute 
$\floor{\log_2 \max\{1,\log |x|\}}$ and the $\ell$ most 
significant bits of $\log x$ in time $O(\ell\log^2\ell)$. \qed 
\end{thm} 

\noindent 
The underlying technique dates back to Gauss and was refined for 
computer use in the 1970s by many researchers 
(see, e.g., \cite{brent,salamin,borwein}).  
We note that in the complexity bound above, we are applying 
the recent $O(n\log n)$ algorithm of Harvey and van der Hoeven for multiplying 
two $n$-bit integers \cite{harvey}. Should we use a more practical (but 
asymptotically slower) integer multiplication algorithm then the time can 
still be kept at $O\!\left(\ell^{1.585}\right)$ or lower. 

Recall that {\em bisection} is the ancient technique of approximating 
a root of a continuous function $f : [r_1,r_2] \longrightarrow 
\R$ by the following trick: If $\sign(f(r_1)f(r_2))\!<\!0$ then 
$f$ must have a root in the open interval $(r_1,r_2)$, and this root lies in 
the left half-interval $\left(r_1,\frac{r_1+r_2}{2}\right)$ if and only 
if $\sign\!\left(f(r_1)f\!\left(\frac{r_1+r_2}{2}\right)\right)\!<\!0$. 
Bisection thus allows one to extract an extra bit of precision for a root of 
$f$ at the cost of one more evaluation of $f$. Put another way, bisection 
allows one to halve the size of an isolating interval at the cost of one more 
evaluation of $f$.  We will need the following result on the bit complexity of 
approximating the real roots of a polynomial in $\Z[x_1]$ by rational numbers. 
\begin{lemma} 
\label{lemma:bisect} 
Suppose $f\!\in\!\Z[x_1]$ has degree $d$, $|f|_\infty\!\leq\!H$, and 
$\ell\!\in\!\N$. Let $\delta(f)$ denote the minimum of $|\zeta_1-\zeta_2|$ 
over all distinct real roots $\zeta_1,\zeta_2$ of $f$. Then, in time \\  
\mbox{}\hfill 
$O\!\left(d^6\log^2(H)+ d^4(d^2+\ell)\ell\log^2(d\ell\log(dH)) \right)$, 
\hfill\mbox{}\\ 
we can find a collection of disjoint non-empty open intervals 
$\{J_i\}^k_{i=1}$ with the following\linebreak properties:\\ 
(a) $k$ is the number of real roots of $f$.\\  
(b) Each $J_i$ contains exactly one root of $f$.\\  
(c) The endpoints of all the $J_i$ are rational numbers 
with logarithmic height\\   
\mbox{}\hspace{.7cm}$O(\ell+d^2+d\log(dH))$.\\ 
(d) All the $J_i$ have width no greater than $2^{-\ell}\delta(f)$. 
\end{lemma} 

\noindent 
{\bf Proof:} The case $\ell\!=\!0$ is well-known in the computational 
algebra community and is elegantly described in \cite{sagraloff}. 
(In fact, \cite{sagraloff} even allows polynomials with real coefficients 
known only up to a given tolerance.) In particular, we merely apply 
Theorem 24 of \cite{sagraloff} to the {\em square-free part}, 
$p\!:=\!f/\gcd(f,f')$, 
of $f$. Note that $p$ can be computed within time $O(d^2\log^2 H)$, 
and its coefficients are integers of logarithmic height 
$\log H'\!=\!O(d+\log H)$, thanks to Lemma \ref{lemma:factor}. So our overall 
complexity bound holds thanks to Lemma \ref{lemma:factor} and the 
$O(d^4\log^2 H')$ bound (in our notation) from \cite[Thm.\ 24]{sagraloff}. 

The case of arbitrary $\ell$ can be derived simply by applying 
bisection after using the $\ell\!=\!0$ case to start with isolating 
intervals that are likely larger than desired, but {\em correct in number}, 
for the real roots of $f$: One first observes that if $\alpha\!\in\!\Q$ has 
logarithmic height $L$ then Proposition \ref{prop:heightsum} implies that 
$f(\alpha)$ has height $O(d\log(H)+d^2L)$. So we can correctly find the sign 
of $f(\alpha)$ by, say, Horner's Rule \cite{vzg}, using 
$O(d\log(H)+d^2L)$ bits of accuracy in all intermediate 
calculations. Since there are at most $d$ 
roots, and each application of Horner's Rule takes $O(d)$ multiplications 
and additions, we see that applying bisection to each of our initial 
isolating intervals (to halve the size of each interval) takes time 
$O((d^3\log(H)+d^4L)\log(dL \log H))$, 
assuming we use the fast multiplication algorithm of \cite{harvey}.  

Corollary \ref{cor:rootspacing} then tells us that 
$|\log\delta(f)|\!=\!O(d^2+d\log(dH))$. This means that 
getting $\ell$ additional bits of accuracy beyond the minimum root 
separation requires time\\ 
\mbox{}\hfill $\sum^{\ell}_{i=0}O\!\left((d^3\log(H)+(L+i)d^4)
\log((L+i)d\log H)\right)$,\hfill\mbox{}\\ 
where $L\!=\!O(d^2+d\log(dH))$. Summing up all our time bounds, we get our 
stated bound. \qed 
\begin{rem} 
We have opted for a streamlined proof at the expense of a 
larger complexity estimate. In particular, the exponent of $d$ 
in our bound can likely be lowered slightly if one uses 
more sophisticted techniques, some of which are discussed 
further in \cite{sagraloff} and the references therein. \dia  
\end{rem}

\section{Our Main Algorithms and Their Complexity} 
\label{sec:algor} 
Our central algorithm for counting roots in $\Rsn$ is conceptually simple but 
ultimately somewhat laborious: Reduce to computing the signs of a linear 
combination of $m$ logarithms at its critical points and poles. To compute 
the signs at the critical points, we simply approximate the input to each 
logarithm, and each resulting summand, to extremely high accuracy. The devil 
is in the level of accuracy, but thanks to our earlier development, the 
required accuracy can be estimated explicitly, and the resulting complexity 
bound is quadratic in $(m\log(BH))^{m+1}$.  
\begin{algor}\label{algor:sign}\mbox{}\\
{\bf Input:} Integers $b_1,\ldots,b_m$ and rational numbers 
$\gamma_{1,1},\gamma_{1,0},\ldots,\gamma_{m,1},\gamma_{m,0},u_0,u_\infty$, 
with $m\!\geq\!2$, $\gamma_{i,1}u+\gamma_{i,0}\!>\!0$ for all 
$u\!\in\!(u_0,u_\infty)$ and $i\!\in\!\{1,
\ldots,m\}$, and $\det\begin{bmatrix}\gamma_{i,1} & 
\gamma_{i,0}\\ \gamma_{j,1} & \gamma_{j,0}\end{bmatrix}\!\neq\!0$ 
for all $i,j$. \\  
{\bf Output:} The signs of $\displaystyle{L(u)\!=\!\sum\limits^m_{i=1}b_i 
\log\left|\gamma_{i,1} u+\gamma_{i,0}\right|}$ at all the critical points of 
$L$ in $(u_0,u_\infty)$.\\  
{\bf Description:}\\
\vspace{-.5cm}
\begin{itemize}
\item[0.]{Let $B\!:=\!\max_i |b_i|$, 
$\log H\!:=\!\max\left\{1,\max_{i,j} h\!\left(\gamma_{i,j}\right)\right\}$,
$\sA := \log^m\left(m^2 16^{m-1}B^2H^{6m-2}\right)$, 
\begin{eqnarray*}
\cE & := & 1.4\cdot m^{6.5} 30^{m+3}(1+\log m)(1+\log(mB))\sA,\\   
\cD & := & m^2e+(m+2)\log\left(8+m2^{m+2}BH^{2m} \right), 
\end{eqnarray*} 
and $\rho\!:=\!1.443(\cD+\log(12m)+\cE)$.} 
\item[1.]{Compute the polynomial 
$\displaystyle{g(u)\!:=\!\sum^m_{i=1}{b_i\gamma_{i,1}\nu_i 
\left. \prod\limits^m_{j=1}(\gamma_{j,1}u +\gamma_{j,0})\nu_j\right/
((\gamma_{i,1}u+\gamma_{i,0})\nu_i)}}$, where $\nu_i$ denotes the least common 
multiple of the denominators of $\gamma_{i,1}$ and $\gamma_{i,0}$.} 
\item[2.]{Via Lemma \ref{lemma:bisect}, find respective isolating intervals 
$J_1,\ldots,J_{k-1}$ to the roots \mbox{$u_1\!<\cdots<\!u_{k-1}$} of 
$g$ in $(u_0,u_\infty)$ such that each $J_i$ has width no greater than 
$2^{-\rho}$. } 
\item[3.]{For all $i\!\in\!\{1,\ldots,k-1\}$ do:} 
\item[4.]{\mbox{}\hspace{.5cm}Let 
$\bar{u}_i\!:=\!\frac{\sup J_i+\inf J_i}{2}$.}
\item[5.]{\mbox{}\hspace{.5cm}For all $j\!\in\!\{1,\ldots,m\}$ do:}
\item[6.]{\mbox{}\hspace{1cm}Compute, via Theorem 
\ref{thm:log}, a rational number $L_j$ agreeing with 
$\log|\gamma_{j,1}\bar{u_i}+\gamma_{j,0}|$\\  
\mbox{}\hspace{1cm}in its first $\ceil{1.443\cE+\log_2(6m)}$ most 
significant bits. } 
\item[7.]{\mbox{}\hspace{.5cm}End For}  
\item[8.]{\mbox{}\hspace{.5cm}Let $\sL_i\!:=\!\sum\limits^m_{j=1}b_jL_j$ and 
$\theta_i\!:=\!\mathrm{sign}(\sL_i)$.}  
\item[9.]{\mbox{}\hspace{.5cm}If $|\sL_i|\!>\!\frac{1}{2}2^{-1.443\cE}$ then} 
\item[10.]{\mbox{}\hspace{1cm}Output ``The sign of $L$ at $u_i$ is 
$\theta_i$.'' } 
\item[11.]{\mbox{}\hspace{.5cm}Else } 
\item[12.]{\mbox{}\hspace{1cm}Output ``$L(u_i)\!=\!0$.'' } 
\item[13.]{\mbox{}\hspace{.5cm}End If } 
\item[14.]{\mbox{}End For}  
\end{itemize}
\end{algor} 
\begin{lemma}
\label{lemma:sign}
Algorithm \ref{algor:sign} is correct and runs in time 
$O(\log(B)+m\log H)^{2m+2}$.   
\end{lemma} 

\noindent 
{\bf Proof:} The correctness of our algorithm follows directly from Theorem 
\ref{thm:bwm} and Lemmata \ref{lemma:crit} and \ref{lemma:spacing}. First 
note that the classical inequality $1-\frac{1}{x}\leq\log x \leq x-1$ (for 
all $x\!>\!0$), yields\\     
\mbox{}\hfill $\frac{s}{v+s} \leq \log(v+s)-\log v \leq \frac{s}{v}$  
\text{ (for all $v\!>\!0$ and $s\!>\!-v$),}  
\hfill\mbox{}\\ 
upon setting $x\!=\!\frac{v+s}{v}$.  
Setting $v\!=\!\gamma_{j,1}\bar{u}_j+\gamma_{j,0}$ and 
$s\!=\!\gamma_{j,1}(u_j-\bar{u}_j)$, and assuming\linebreak  
$\gamma_{j,1}u_j+\gamma_{j,0},\gamma_{j,1}\bar{u}_j+\gamma_{j,0}\!>\!0$, 
we then obtain 
\begin{eqnarray}
\label{ineq:trick} 
\frac{(u_j-\bar{u}_j)\gamma_{j,1}}{\gamma_{j,1}u_j+\gamma_{j,0}}
\leq \log(\gamma_{j,1}u_j+\gamma_{j,0})-
\log(\gamma_{j,1}\bar{u}_j+\gamma_{j,0})
\leq \frac{(u_j-\bar{u}_j)\gamma_{j,1}}{\gamma_{j,1}\bar{u}_j+\gamma_{j,0}}.  
\end{eqnarray} 
The proof of Assertion (3) of Lemma \ref{lemma:spacing} tells us that 
$\frac{|\gamma_{j,1}u_j+\gamma_{j,0}|}{|\gamma_{j,1}|}\!\geq\!\Delta$. 
Since $1/\log 2\!<\!1.443$ we have $\Delta\!>\!2^{-1.443\cD}$, thanks  
to the definition of $\cD$. So the 
definition of $s$ tells us that $|u_j-\bar{u}_j|\!\leq\!\frac{1}{2}
2^{-1.443\cD}$ is sufficient to guarantee that 
$\frac{|\gamma_{j,1}\bar{u}_j+\gamma_{j,0}|}{|\gamma_{j,1}|}\!\geq\!\Delta/2$. 
So, by Inequality (\ref{ineq:trick}), we obtain that 
$|u_j-\bar{u}_j|\!\leq\!2^{-\rho}$ guarantees 
$|\log(\gamma_{i,1}u_i+\gamma_{i,0})-\log(\gamma_{i,1}\bar{u}_i+\gamma_{i,0})|
\!\leq\!\frac{1}{6m}2^{-1.443\cE}$. Should $\gamma_{i,1}u+\gamma_{i,0}\!<\!0$ 
we can repeat our preceding argument, with a sign flip, to obtain 
that $|u_j-\bar{u}_j|\!\leq\!2^{-\rho}$ guarantees 
$|\log|\gamma_{i,1}u_i+\gamma_{i,0}|-\log|\gamma_{i,1}\bar{u}_i+\gamma_{i,0}||
\!\leq\!\frac{1}{6m}2^{-1.443\cE}$. 
So then, thanks to Step 6 and the Triangle Inequality, we see that our 
algorithm computes, for each $i\!\in\!\{1,\ldots,k-1\}$, a rational $\sL_i$ 
such that $|L(u_i)-\sL_i|\!\leq\!\frac{1}{3}2^{-1.443\cE}$. 

Theorem \ref{thm:bwm} then tells us that $|L(u_i)|$ is either $0$ or strictly 
greater than $2^{-1.443\cE}$. So the threshold on $|\sL_i|$ from Step 9 indeed 
correctly distinguishes between $L(u_i)$ being nonzero or zero, and the signs 
of $L(u_i)$ and $\sL_i$ also match when $L(u_i)\!\neq\!0$ thanks to our chosen 
accuracy. In other words, our algorithm is correct.   

We now analyze the complexity of our algorithm. 
First note that $H$, $\sA$, $\cE$, $\cD$, and $\rho$ need not be 
computed exactly: it is sufficient to work with the ceilings of these 
quantities, or even the smallest powers of $2$ respectively greater than these 
quantities. In particular, these parameters can easily be computed via 
standard efficient methods for computing the exponential function 
\cite{ahrendt} (along with 
Theorem \ref{thm:log}) and thus the complexity of Step 0 is negligible, and 
in fact asymptotically dominated by Steps 2 and beyond. 

Likewise, Step 1 is easily seen to take time within 
$O(m^2(\log^2(B)+m\log^2 H))$, to be asymptotically dominated by Step 2 and 
beyond. (The preceding bound is a crude over-estimate of what can be obtained 
by combining the fast polynomial multiplication method from, say, 
\cite[Sec.\ 8.4]{vzg} with the fast integer multiplication method of Harvey 
and van der Hoeven \cite{harvey}.)  

Lemma \ref{lemma:crit} tells us that the complexity of Step 2 can be estimated 
by replacing $(d,H,\ell)$ in the statement of Lemma \ref{lemma:bisect} by 
$(m-1,m2^{m-1}BH^{2m},\rho)$. Noting that $\rho\!=\!O(\cE)$ and 
$m^r,\log^r B, \log^r H\!=\!O(\cE)$ for any $r\!>\!0$, Lemma \ref{lemma:bisect} 
then tells us that Step 2 takes time $O(m^4\cE^2\log^2 \cE)$, which is 
bounded from above by
\begin{eqnarray} 
\label{eqn:biggest} 
& O\!\left(m^{17}\log^2(m)900^m\log^2(mB)\right)  
\cdot \left(O(\log(B)+m\log H)\right)^{2m} & \nonumber \\ 
& = O\!\left(901^m\log^2(mB)\right)
\cdot \left(O(\log(B)+m\log H)\right)^{2m} & \nonumber \\ 
& = O(\log(B)+m\log H)^{2m+2} & 
\end{eqnarray} 

Thanks to Theorem \ref{thm:log}, a simple over-estimate for the 
complexity of Step 6 is 
\begin{eqnarray*} 
O(m^2\log^2(mB)\log^2(m\log(BH)))\cdot O(\log(B)+m\log H)^m, 
\end{eqnarray*}   
so then the time spent in (each run of) Steps 4--7 in total is 
\begin{eqnarray*} 
O(m^3\log^2(mB)\log^2(m\log(BH)))\cdot O(\log(B)+m\log H)^m.  
\end{eqnarray*}   

Since $k-1\!\leq\!m-1$, Steps 3--14 then take time no greater than 
\begin{eqnarray*} 
O(m^4\log^2(mB)\log^2(m\log(BH)))\cdot O(\log(B)+m\log H)^m.  
\end{eqnarray*}   

We thus see, from comparison to Estimate (\ref{eqn:biggest}), 
that Step 2 in fact dominates the asymptotic complexity of our 
entire algorithm, so our complexity estimate follows. \qed 

\medskip 
We can now state our algorithm for counting the positive roots of circuit 
systems: 
\begin{algor}\label{algor:main}\mbox{}\\
{\bf Input:} \scalebox{.95}[1]{Polynomials $f_1,\ldots,
f_n\!\in\!\Z\!\left[x^{\pm 1}_1,\ldots, x^{\pm 1}_n\right]$ with 
$A\!:=\!\bigcup_i \supp(f_i)$ a circuit and $\#A\!=\!n+2$.}\\  
{\bf Output:} The number of roots of $F\!=\!(f_1,\ldots,f_n)$ in $\Rn_+$.\\  
{\bf Description:}\\
\vspace{-.5cm}
\begin{itemize}
\item[0.]{ Find the unique non-degenerate sub-circuit $\Sigma$ of $A$,  
re-index the points of $A$ so that $\Sigma\!=\!\{a_1,\ldots,a_m,a_{n+1},
a_{n+2}\}$, translate $a_1,\ldots,a_{n+1}$ by $-a_{n+2}$, and 
set $a_{n+2}\!:=\!\bO$.}  
\item[1.]{ Letting $[c_{i,j}]$ be the coefficient matrix of $F$, 
check whether all the $n\times n$ sub-matrices of
\scalebox{.7}[.7]{$\begin{bmatrix}c_{1,1} & \cdots & c_{1,n} & c_{1,n+2}\\ 
\vdots & \ddots & \vdots \\ 
c_{n,1} & \cdots & c_{n,n} & c_{n,n+2}\end{bmatrix}$} are non-singular, 
and whether \scalebox{.7}[.7]{$\begin{bmatrix}c_{i,n+1} & c_{i,n+2}\\ 
c_{j,n+1} & c_{j,n+2}\end{bmatrix}$} is non-singular for all $i\!\neq\!j$. 
If not, then Output\\ 
``Your system might have infinitely many roots but I'm not sure: Please 
check if there are any updates to this algorithm, addressing the cases of  
vanishing minors.''\\  
and {\tt STOP}. } 
\item[2.]{Reduce $F\!=\!\bO$ to a system of equations of the form
$G\!=\!\bO$, where $G\!:=\!(g_1,\ldots,g_n)$ and 
$g_i(x)\!:=\!x^{a_i}-\gamma_{i,1}x^{a_{n+1}}-\gamma_{i,0}$ for all $i$. }
\item[3.]{ Let $b\!\in\!\Z^{(m+2)\times 1}$ be the unique (up to sign)
minimal circuit relation of $\Sigma$, and define\\ 
$\displaystyle{L(u)\!:=\!b_{m+1}\log|u|
+\sum\limits^m_{i=1} b_i\log|\gamma_{i,1}u+\gamma_{i,0}|}$ and\\ 
\mbox{}\hfill  
$\displaystyle{I\!:=\!\{u\!\in\!\R_+ \; | \; \gamma_{i,1}u+\gamma_{i,0}\!>\!0 
\text{ for all } i\!\in\!\{1,\ldots,n\}\}}$.\hfill\mbox{} } 
\item[4.]{ Via Algorithm \ref{algor:sign} and Proposition \ref{prop:crit},  
compute the number, $N$, of roots of $L$ in $I$, and Output $N$.}     
\end{itemize} 
\end{algor} 

\begin{lemma}
\label{lemma:main}
Algorithm \ref{algor:main} is correct and runs in time 
$O(n\log(nd)+n^2\log H)^{2n+4}$, where $d$ is the largest absolute 
value of an entry of $A$ and $H\!:=\!\max_i|f_i|_\infty$.    
\end{lemma} 

\noindent 
{\bf Proof:} Letting $u_0\!:=\!\inf I$, $u_\infty\!:=\!u_k\!:=\!\sup I$, 
and letting $u_1\!<\cdots<\!u_{k-1}$ be the critical points of $L$ in 
$I$ as before, note that the sign of 
$\lim_{u\rightarrow u^+_0} L(u)$ is either $-b_{m+1}$ (if $u_0\!=\!0$) 
or $-b_{i'}$ with $i'$ the unique index such that 
$u_0\!=\!-\gamma_{i',0}/\gamma_{i',1}$. (This index is unique thanks 
to Step 1 ensuring the $2\times 2$ minor condition.) Similarly,  
$\lim_{u\rightarrow u^-_k} L(u)$ is simply the sign of $-b_{i''}$ 
(resp.\ $b_1+\cdots+b_{m+1}$) for some easily computable unique index 
$i''\!\in\!\{1,\ldots,m+1\}$ if $u_k\!<\!\infty$ (resp.\ $u_k\!=\!+\infty$). 
So the use of Proposition \ref{prop:crit} is clear.  

The correctness of Algorithm \ref{algor:main} then follows 
directly from Lemmata \ref{lemma:gale}, Proposition \ref{prop:crit}, 
and Lemma \ref{lemma:sign}. So we now analyze the complexity of our 
algorithm. 

Thanks to Lemmata \ref{lemma:lin} and \ref{lemma:rightnull}, 
it is clear that Steps 0--3 are doable in time\linebreak  
$O(n^{2+\omega}\log^2(ndH))$. This will not be the dominant 
part of the algorithm: Observing\linebreak 
that $h(\gamma_{i,j})\!=\!O(n\log(nH))$ and $h(b_i)\!=\!O(n\log(nd))$ for all 
$i,j$, Lemma \ref{lemma:sign} then tells\linebreak 
us that applying Algorithm \ref{algor:sign} and Proposition 
\ref{prop:crit} (with $m\!\leq\!n+1$) takes time\linebreak 
$O(n\log(nd)+n^2\log(nH))^{2n+4}$. \qed 

\medskip 
We are now ready to state the analogue of Lemma \ref{lemma:main} 
for counting roots in $\Rsn$: 
\begin{lemma}
\label{lemma:main+} 
Given any $(n+2)$-nomial $n\times n$ system $F\!=\!(f_1,\ldots,
f_n)\!\in\!\Z\!\left[x^{\pm 1}_1,\ldots, x^{\pm 1}_n\right]^n$ supported on 
a circuit $A$ with cardinality $n+2$, we can count exactly the number of roots 
of $F$ in $\Rsn$ in time $O(n\log(nd)+n^2\log H)^{2n+4}$, where $d$ is the 
largest absolute value of an entry of $A$ and $H\!:=\!\max_i|f_i|_\infty$.    
\end{lemma} 

\noindent 
{\bf Proof:} The proof is almost identical to that of Lemma \ref{lemma:main}, 
save that we apply Corollary \ref{cor:crit} instead of Proposition 
\ref{prop:crit}, Lemma \ref{lemma:gale+} instead of Lemma \ref{lemma:gale}, 
and that we use a {\em modified} version of Algorithm \ref{algor:main}. 

In particular, the modifications to Algorithm \ref{algor:main} are that 
we replace Steps 0, 3, and 4 respectively by Steps 0', 3', and 4' stated below:\\  
\fbox{\vbox{
\noindent 
0'. {\em Find the unique non-degenerate sub-circuit $\Sigma$ of $A$ and 
underlying minimal circuit relation\\
\mbox{}\hspace{.6cm}$b\!\in\!\Z^{(m+2)\times 1}$, 
re-index the points of $A$ so that $\Sigma\!=\!\{a_1,\ldots,a_m,a_{n+1},
a_{n+2}\}$ and $b_{m+1}$ is\\ 
\mbox{}\hspace{.6cm}odd, translate $a_1,\ldots,a_{n+1}$ by $-a_{n+2}$, and
set $a_{n+2}\!:=\!\bO$.}}} 

\noindent 
\fbox{\mbox{3'. {\em Define $\displaystyle{L(u)\!:=\!b_{m+1}\log|u|
+\sum\limits^m_{i=1} b_i\log|\gamma_{i,1}u+\gamma_{i,0}|}$}}}
  
\noindent 
\fbox{\vbox{
\noindent 
4'. {\em Via Algorithm \ref{algor:sign}, compute the number 
$\cN$ from Corollary \ref{cor:crit} and the number of\\  
\mbox{}\hspace{.5cm}
degenerate roots of $L$ in $\R$. Output their sum.}}}

\smallskip 
\noindent 
Note in particular that, in Step 0', $b$ must have an odd coordinate since 
minimal circuit relations are assumed to have relatively prime coordinates. 
Also, the left or right-handed limits of $L$ at a real (possibly infinite) 
pole are easy to compute via the sign of a suitable $b_i$ (if the pole is 
finite) or the sign of $b_1+\cdots+b_{m+1}$ (if the pole is $\pm \infty$). 
The correctness of the modified algorithm is then immediate. 

The complexity analysis for our modified algorithm is almost identical to that 
of Algorithm \ref{algor:main}, save that there is extra work taking time 
$O(n\cdot n^2)$ to compute the signs of $\Lambda(u)$ and the $\Gamma'(u_i)$. 
This is negligible compared to the other steps, so our final asymptotic 
complexity bound remains the same. \qed 

\begin{rem} 
While we could have simply applied Algorithm \ref{algor:main} $2^n$ times 
(once for each orthant of $\Rsn$), our proof above is 
more practical: It enables complexity polynomial in $n$ for 
root counting in $\Rsn$ should 
sufficiently sharp new diophantine estimates become available 
(see Remark \ref{rem:abc} from the introduction). \dia 
\end{rem} 

\section{Affine Roots and Proving Theorem \ref{thm:main}}  
\label{sec:proof} 
To conclude, we'll first need to recall some observations regarding sparse 
resultants and roots of $t$-nomial $(n+1)\times n$ (over-determined) systems 
on coordinate subspaces. For further details on sparse resultants we refer the 
reader to \cite{gkz94}. 
\begin{prop} 
\label{prop:intersect} 
If $A\!\subset\!\Zn$ is a circuit and $X$ is any coordinate 
subspace of $\Rn$ then $A\cap X$ is either empty, the vertex 
set of a simplex, or a circuit. \qed 
\end{prop} 
\begin{lemma} 
\label{lemma:res} 
\cite[Ch.\ 8]{gkz94} 
Suppose $A\!=\!\{a_1,\ldots,a_t\}\!\subset\!\Zn$ has cardinality $t\!\leq\!n+2$ 
and does not lie in any affine hyperplane, and $F\!=\!(f_1,\ldots,f_{n+1})$ 
with $f_i(x)\!=\!\sum^{n+1}_{j=1} c_{i,j} x^{a_j}$ for all $i$ and 
the $c_{i,j}$ algebraically independent indeterminates. 
Then there is a unique (up to sign) polynomial\\ 
\mbox{}\hfill 
$\Delta_A\!\in\!\Z\!\left[c_{i,j}\; | \; (i,j)\!\in\!\{1,\ldots,n\}\times\{1,
\ldots,t\}\right]\setminus\{0\}$,\hfill\mbox{}\\ 
of minimal degree, such that $\left[[c_{i,j}]\!\in\!\C^{n\times t}\right.$ and 
$\left.\Delta_A(\ldots,c_{i,j},\ldots)\!\neq\!0\right] \Longrightarrow F$ has 
no roots in $\Csn$. In particular, if $t\!=\!n+1$ (resp.\ $t\!=\!n+2$) then 
$\Delta_A\!=\!\det[c_{i,j}]$ (resp.\ $\deg \Delta_A\!\leq\!n\cdot n!
V$),\linebreak 
\scalebox{.92}[1]{where $V$ is the volume of the convex hull of $A$, 
normalized so that the unit $n$-cube has volume $1$. \qed}  
\end{lemma} 

The polynomial $\Delta_A$ above is an example of a {\em sparse resultant}, and 
is one of many ways to formulate the fact that $(n+1)$-tuples of $n$-variate 
polynomials generically have no roots. Our statement above is minimalist 
and is meant to approach affine roots in the following simple way: 
\begin{lemma} 
\label{lemma:affine} 
Suppose 
$F\!=\!(f_1,\ldots,f_n)\!\in\!\left(\C\!\left[x^{\pm 1}_1,\ldots,x^{\pm 1}_n
\right] \right)^n$ is a $t$-nomial $n\times n$ polynomial system supported on 
a fixed $A\!\subset\!\Zn$ and $\supp(f_i)\!=\!A$ for all $i$.  
Let $\pi_i : \Rn \longrightarrow \R$ denote the projection sending 
$(x_1,\ldots,x_n)\mapsto x_i$. Then:\\ 
1.\ For any $i\!\in\!\{1,\ldots,n\}$ we have that 
$\min\pi_i(A)\!>\!0 \Longrightarrow F$ vanishes on all of\\ 
\mbox{}\hspace{.5cm}$(\Cs)^{i-1}\times\{0\} \times(\Cs)^{n-i}$.\\  
2.\ If $\min\pi_i(A)\!\geq\!0$ for all $i\!\in\!\{1,\ldots,\ell\}$ 
and no points of $A$ lie on the coordinate subspace\\ 
\mbox{}\hspace{.5cm}$\{x_1\!=\cdots=\!
x_\ell\!=\!0\}$ then $F$ vanishes on all of $\{0\}^\ell\times 
(\Cs)^{n-\ell}$.\\ 
3.\ For any $i\!\in\!\{1,\ldots,n\}$ we have that $\min\pi_i(A)\!<\!0 
\Longrightarrow F$ has {\em no} complex roots in the\\ 
\mbox{}\hspace{.5cm}coordinate hyperplane $\{x_i\!=\!0\}$.\\ 
4.\ If $\min\pi_i(A)\!\geq\!0$ for all $i\!\in\!\{1,\ldots,\ell\}$ 
and $A\cap\{x_1\!=\!\cdots=\!x_\ell\!=\!0\}\!\neq\!\emptyset$  
then $F$ generically has\linebreak 
\mbox{}\hspace{.5cm}{\em no} roots in $\{0\}^\ell\times (\Cs)^{n-\ell}$. 
\end{lemma} 
\begin{ex} 
It is easy to check (by specializing subsets of variables to $0$, 
and then dividing out by $x_1x_2x_3$)  
that the real zero set of the $4$-nomial $3\times 3$ system 
\begin{eqnarray*} 
\left( x_1x_2 + x_2x_3  +  x_3x_1 - 3x_1x_2x_3, \right. \\  
   \mbox{}x_1x_2 + 2x_2x_3 + 4x_3x_1 - 7x_1x_2x_3,  \\  
\left. x_1x_2 + 3x_2x_3 + 9x_3x_1 - 13x_1x_2x_3 \right)  
\end{eqnarray*} 
is the union of the point $(1,1,1)$, the $x_1$-axis, the 
$x_2$-axis, and the $x_3$-axis. Note in particular that while 
the underlying support $A$ satisfies $\min\pi_i(A)\!=\!0$ for all $i$,  
$A$ does not intersect any of the coordinate axes. Hence the 
infinitude of real roots for our example. \dia 
\end{ex} 

\noindent 
{\bf Proof of Lemma \ref{lemma:affine}:} Assertion (1) is merely the case of 
all the $f_j$ being divisible by $x_i$. Assertion (2) is immediate upon 
rewriting the $f_j$ as elements of 
$\left(\C\!\left[x^{\pm 1}_{\ell+1},\ldots,x^{\pm 1}_n\right]\right)
[x_1,\ldots,x_\ell]$: The condition on $A$ holds if and only if 
all the $f_j$ lack a non-trivial monomial term of degree $0$ with respect to 
$(x_1,\ldots,x_\ell)$. Assertion (3) is merely the observation that 
roots with $x_i\!=\!0$ are impossible if some $f_j$ has a monomial 
with negative exponent for $x_i$. 

To prove Assertion (4) observe that setting $x_1\!=\cdots=\!x_\ell\!=\!0$ 
in $F$ results in such an $F$ generically becoming an (over-determined) 
$n\times (n-\ell)$ 
system. In particular, for $F$ to fail to have roots, it is 
enough for the first $n-\ell+1$ polynomials to fail to have a common 
root. So by Lemma \ref{lemma:res} we are done. \qed 

\subsection{The Proof of Theorem \ref{thm:main}} 
\label{sub:proof} 
The cases of root counting in $\Rn_+$ and $\Rsn$ follow respectively from 
Lemmata \ref{lemma:main} and \ref{lemma:main+} when $t\!=\!n+2$. For 
$t\!=\!n+1$ we simply use Lemma \ref{lemma:n+1} and Corollary \ref{cor:n+1+} 
instead. 

To count roots in $\Rn$, let us first compute $\min \pi_i(A)$ for all 
$i\!\in\!\{1,\ldots,n\}$: This can be done simply by reading the coordinates 
of the points of $A$, and this takes time $O(n^2\log d)$. Let $A_i$ denote the 
intersection of $A$ with the $x_i$-axis. Then we can also easily compute 
$A_1,\ldots,A_n$ and $A_1\cap \cdots \cap A_n$ within the same time bound as 
well. Since we have already derived how to efficiently count roots in $\Rsn$ 
and $\Rn_+$, it suffices to prove that we can count roots on all the 
proper coordinate subspaces of $\Rn$ within our stated time bound. 

By Assertion (3) of Lemma \ref{lemma:affine}, $\min\pi_i(A)\!<\!0$ for all $i$ 
implies that all the roots of $F$ in $\Rn$ lie in $\Rsn$. So, without loss of 
generality, we can re-order variables and assume $\min\pi_i(A)\!<\!0$ for all 
$i\!\geq\!\ell+1\!\geq\!2$. 

By Assertion (1) of Lemma \ref{lemma:affine}, if $\min\pi_i(A)\!>\!0$ for some 
$i\!\in\!\{1,\ldots,\ell\}$ then $F$ has infinitely many real roots. So we 
may assume that $\min\pi_i(A)\!=\!0$ for all $i\!\in\!\{1,\ldots,\ell\}$. 

Observe now that if $\ell\!=\!n$ (i.e., if $\min\pi_i(A)\!=\!0$ for all $i$) 
and $A_i$ is empty for some $i$ then $F$ has infinitely many real roots by 
Assertion (2) of Lemma \ref{lemma:affine}. So if $\ell\!=\!n$ we may assume 
that $A_i$ is non-empty for all $i$. Since $\supp(f_i)\!=\!A$ generically for 
all $i$, Assertion (4) of Lemma \ref{lemma:affine} tells us that $F$ 
generically has no 
roots on any coordinate subspace containing one of the coordinate axes. So 
now we are left with checking whether $F$ vanishes at $\bO$. Since $F$ 
vanishes at $\bO$ if and only if $A_1\cap \cdots \cap A_n$ is non-empty, 
the case $\ell\!=\!n$ is done and we assume now that $\ell\!<\!n$. In other 
words, we now need only check for roots on the coordinate subspace 
$\{x_1\!=\cdots=\!x_\ell\!=\!0\}$ of $\Rn$. 

Now, should $A_1\cap\cdots\cap A_\ell$ be empty, then Assertion (2) 
of Lemma \ref{lemma:affine} tells us that $F$ has infinitely many real roots. 
So we may assume $A_1\cap\cdots \cap A_\ell\!\neq\!\emptyset$. Since 
$\supp(f_i)\!=\!A$ generically, Assertion (4) of Lemma \ref{lemma:affine} then 
tells us that $F$ generically has no roots on any coordinate subspace of $\Rn$ 
containing $\{x_1\!=\cdots=\!x_\ell\!=\!0\}$. 

We have thus shown how to count roots in all coordinate subspaces of $\Rn$. 
In particular, we see that the complexity of counting roots in $\Rn$ is 
asymptotically the same as our complexity bound for counting in $\Rsn$, since 
the only additional work is computing the $\pi_i(A)$, the $A_i$, 
$A_1\cap\cdots \cap A_\ell$, and $A_1\cap\cdots\cap A_n$, which takes 
time neglible compared to our main bound. \qed 

\medskip 
A consequence of our preceding proof is that the additional genericity 
condition needed to efficiently count roots in $\Rn$ can be taken to be 
the following:\\ 
\mbox{}\hspace{.5cm}{\em All the coefficients of $F$ are nonzero, and a 
collection of sparse resultants (one for each\\
\mbox{}\hspace{.5cm}coordinate subspace of $\Rn$), in the coefficients of $F$, 
all not vanish.}\\ 
So by Proposition \ref{prop:intersect} and Lemma \ref{lemma:res}, 
the degree of the product of all these underlying polynomials is no greater  
than\\ 
\mbox{}\hfill 
$n(n+2)+d^n+nd^{n-1}+\binom{n}{2}d^{n-2}+\cdots+nd\!=\!n^2+n+(d+1)^n-1$. 
\hfill\mbox{}\\  
By Schwartz's Lemma \cite{schwartz}, as applied in the proof of Corollary 
\ref{cor:gen}, it is then clear that $F$ having uniformly random coefficients 
in $\left\{-\ceil{\frac{n^2+n+(d+1)^n-1}{2\eps}},\ldots, 
\ceil{\frac{n^2+n+(d+1)^n-1}{2\eps}}\right\}$ is enough to make our 
additional genericity condition hold with probability greater than $1-\eps$.

\section*{Acknowledgements} 
I thank Dan Bates and Jon Hauenstein for answering my questions on  
how {\tt Bertini} handles polynomial systems of extremely high degree. I also 
thank Jan Verschelde for answering my questions on fine-tuning the options in 
{\tt PHCpack}. Special thanks to Timo de Wolff for pointing out reference 
\cite{chandra} and Alexander Barvinok for pointing out reference \cite{proy}. 

I am also happy to thank Weixun Deng, Alperen Erg\"{u}r, and Grigoris 
Paouris for good company and inspirational conversations.  

\section*{Farewell to a Friend} 
Tien-Yien Li passed away a few months into the COVID-19 pandemic. 
TY (as he was known to his friends) was an immensely kind and generous 
man, and a dear friend, in addition to being a great mathematician. 
Through hours-long grilling sessions in October 1993, at the Centre de Recerca 
Matematica in Barcelona, he taught me lessons on perseverance, curiosity, 
scholarship, and generosity that I would always remember. It was there that 
I also got to know TY and his unique sense of humor. He always faced the 
greatest difficulties with a smile. 
I admired him both as a person and a mathematician. I truly miss him.  

\bibliographystyle{amsalpha}

\begin{thebibliography}{A}    

\bibitem[Ahr99]{ahrendt} Timm Ahrendt, {\it ``Fast computations of the
exponential function,''} in proceedings of STACS '99 (16th
annual conference on Theoretical aspects of computer science),
pp.\ 302--312, Springer-Verlag Berlin, 1999.  


\bibitem[BS96]{bachshallit} Eric Bach and Jeff Shallit, {\it
Algorithmic Number Theory, Vol.\ I: Efficient Algorithms,}
MIT Press, Cambridge, MA, 1996.

\bibitem[Bak77]{baker} Alan Baker, {\it ``The Theory of Linear 
Forms in Logarithms,''} in Transcendence Theory: Advances and Applications:
proceedings of a conference held at the University of Cambridge,
Cambridge, Jan.--Feb., 1976, Academic Press, London, 1977. 

\bibitem[Bak98]{bakerabc} Alan Baker, {\it ``Logarithmic forms and the
$abc$-conjecture,''} Number theory (Eger, 1996), pp.\ 37--44, de Gruyter,
Berlin, 1998. 

\bibitem[BW93]{bakerwustholtz} Alan Baker and Gisbert Wustholtz, 
{\it ``Logarithmic forms and group varieties,''} J.\ Reine Angew.\ 
Math.\ {\bf 442} (1993), pp.\ 19--62.

\bibitem[BPR06a]{bprdim} Saugata Basu; Richard Pollack; 
Marie-Fran\c{c}oise Roy, {\it ``Computing the dimension of a semi-algebraic 
set,''} reprinted in J.\ Math.\ Sci.\ (N.Y.) 134 (2006), no.\ 5, pp.\ 
2346--2353.   

\bibitem[BPR06b]{bpr} Saugata Basu; Richard Pollack; and Marie-Fran\c{c}oise 
Roy, {\it Algorithms in real algebraic geometry,} 2nd edition, 
Algorithms and Computation in Mathematics, 10, Springer-Verlag, Berlin, 2006

\bibitem[BHNS16]{batesgale} Daniel J.\ Bates; Jonathan D.\ Hauenstein; 
Matthew E.\ Niemerg; Frank Sottile, {\it ``Software for the Gale transform of 
fewnomial systems and a Descartes rule for fewnomials,''} 
Numer.\ Algorithms 73 (2016), no.\ 1, pp.\ 281--304.  

\bibitem[BHSW13]{bertini} Daniel J.\ Bates, Jon D.\ Hauenstein, Andrew J.\ 
Sommese, and Charles W.\ Wampler, {\it Numerically solving polynomial systems 
with Bertini,} Software, Environments, and Tools 25, SIAM, 2013.  


\bibitem[BKR86]{bkr} Michael Ben-Or; Dexter Kozen; and John Reif, {\it 
``The Complexity of Elementary Algebra and Geometry,''} J.\ Computer and 
System Sciences 32 (1986), pp.\ 251--264.  

\bibitem[Ber03]{dan} Daniel J.\ Bernstein, {\it ``Computing
Logarithm Intervals with the Arithmetic-Geometric Mean Iterations,''}
available from {\tt http://cr.yp.to/papers.html } .

\bibitem[BBS06]{bbs} Benoit Bertrand; Fr\'{e}d\'{e}ric Bihan; and
Frank Sottile, {\it ``Polynomial Systems with Few Real
Zeroes,''} Mathematisches Zeitschrift, 253 (2006), no.\ 2, pp.\ 361--385.

\bibitem[Bih11]{bihan} Fr\'{e}d\'{e}ric Bihan, {\it Topologie des 
vari\'{e}t\'{e}s creuses,} Habilitation thesis, Universit\'{e} de 
Savoie, France, 2011. 

\bibitem[BD17]{bihandickenstein} Fr\'{e}d\'{e}ric Bihan and Alicia 
Dickenstein, {\it ``Descartes’ Rule of Signs for Polynomial Systems Supported 
on Circuits,''} International Mathematics Research Notices, Vol.\ 2017, Issue 
22, November 2017, Pages 6867–6893.  

\bibitem[BDF20]{optcktdescartes} Fr\'{e}d\'{e}ric Bihan; Alicia Dickenstein; 
Jens Forsg\aa{}rd, {\it ``Optimal Descartes' rule of signs for systems 
supported on circuits,''} Math ArXiV preprint {\tt arXiv:2010.09165v1 } . 

\bibitem[BDG19]{sparsebio1} Fr\'{e}d\'{e}ric Bihan; Alicia Dickenstein; 
and Magal\'{\i} Giaroli, {\it ``Regions of multistationarity in cascades of 
Goldbeter-Koshland loops,''}  
J.\ Math.\ Biol.\ (2019) Vol.\ 78(4), pp.\ 1115--1145.

\bibitem[BDG20a]{sparsechem2} Fr\'{e}d\'{e}ric Bihan; Alicia Dickenstein;
and Magal\'{\i} Giaroli, {\it ``Lower bounds for positive roots and regions of 
multistationarity in chemical reaction networks,''} J.\ Algebra (2020), Vol.\ 
542, pp.\ 367-411.

\bibitem[BDG20b]{signone} Fr\'{e}d\'{e}ric Bihan; Alicia Dickenstein; and 
Magal\'{\i} Giaroli, {\it ``Sign conditions for the existence of at least one 
positive solution of a sparse polynomial system,''} Advances in Mathematics, 
to appear, 2020. 

\bibitem[BRS09]{brs} Fr\'{e}d\'{e}ric Bihan; J.\ Maurice Rojas; Casey E.\ 
Stella, {\it ``Faster Real Feasibility via Circuit Discriminants,''}
proceedings of International Symposium on Symbolic and Algebraic Computation
(ISSAC 2009, July 28--31, Seoul, Korea), pp.\ 39--46, ACM Press, 2009. 

\bibitem[BS07]{bs} Fr\'{e}d\'{e}ric Bihan and Frank Sottile, {\it ``New
Fewnomial Upper Bounds from Gale Dual Polynomial Systems,''}
Moscow Mathematical Journal, 7 (2007), no.\ 3, pp.\ 387--407.

\bibitem[BCSS98]{bcss} Lenore Blum; Felipe Cucker; Mike Shub; and
Steve Smale, {\it Complexity and Real Computation,} Springer-Verlag, 1998.  

\bibitem[BG06]{bombgu} Enrico Bombieri and Walter Gubler, {\it Heights 
in Diophantine Geometry,} new mathematical monographs: 4, Cambridge University 
Press, 2006. 

\bibitem[BBK09]{bomb} Enrico Bombieri; Jean Bourgain; and Sergei Konyagin,
{\it ``Roots of polynomials in subgroups of $\F^*_p$ and
applications to congruences,''} Int.\ Math.\ Res.\ Not.\ IMRN 2009,
no.\ 5, pp.\ 802--834. 

\bibitem[BB88]{borwein} John M.\ Borwein and Peter B.\ Borwein; {\it ``On
the Complexity of Familiar Functions and Numbers,''}
SIAM Review, Vol.\ 30, No.\ 4, (Dec., 1988), pp.\ 589--601.

\bibitem[Bre76]{brent} Richard P.\ Brent, {\it ``Fast Multiple-Precision
Evaluation of Elementary Functions,''} Journal of the Association
for Computing Machinery, vol.\ 23, No.\ 2, April 1976, pp.\ 242--251.

\bibitem[BMS06]{bms} Yan Bugeaud; Maurice Mignotte; and Samir Siksek, 
{\it ``Classical and modular approaches to exponential Diophantine 
equations, I, Fibonacci and Lucas perfect powers,''}  Ann.\ of Math.\ 
(2) {\bf 163} (2006), pp.\ 969--1018. 

\bibitem[BET19]{bet} Peter B\"{u}rgisser, Alperen A.\ Erg\"{u}r, and 
Josu\'{e} Tonelli-Cueto, {\it ``On the Number of Real Zeros of Random 
Fewnomials,''} SIAM Journal on Applied Algebra and Geometry, 3(4), 
pp.\ 721--732, 2019. 

\bibitem[CFKLLS00]{cfklls} Ran Canetti; John B.\ Friedlander;
Sergey Konyagin; Michael Larsen; Daniel Lieman; and Igor E.\ Shparlinski,
{\it ``On the statistical properties of Diffie-Hellman distributions,''}
Israel J.\ Math.\ 120 (2000), pp.\ 23--46. 

\bibitem[Can88]{canny} John F.\ Canny, {\it ``Some Algebraic
and Geometric Computations in PSPACE,''} Proc.\ 20$\thth$ ACM
Symp.\ Theory of Computing, Chicago (1988), ACM Press.   

\bibitem[CD07]{catdick} Eduardo Cattani and Alicia Dickenstein, {\it 
``Counting solutions to binomial complete intersections,''} Journal 
of Complexity 23 (2007), pp.\ 82--107. 

\bibitem[CS16]{chandra} Venkat Chandrasekaran and Parikshit Shah, 
{\it ``Relative Entropy Relaxations for Signomial Optimization,''} 
SIAM J.\ Optim., Vol.\ 26, No.\ 2, pp.\ 1147--1173, 2016.  

\bibitem[CL14]{tianran} Tianran Chen and Tien-Yien Li, {\it ``Solutions
to Systems of Binomial Equations,''}
Annales Mathematicae Silesianae {\bf 28} (2014), pp.\ 7--34.

\bibitem[CGRW17]{cgrw} Qi Cheng; Shuhong Gao; J.\ Maurice Rojas; and
Daqing Wan, {\it ``Sparse Univariate Polynomials with Many Roots Over a
Finite Field,''} Finite Fields and their Applications, Vol.\ 46, July 2017,
pp.\ 235--246. 

\bibitem[CG84]{chigo} Alexander L.\ Chistov and Dima Yu Grigoriev, {\it
``Complexity of Quantifier Elimination in the Theory of Algebraically
Closed Fields,''} Lect.\ Notes Comp.\ Sci.\ 176, Springer-Verlag (1984).

\bibitem[CLRS09]{clrs} Thomas H.\ Cormen; Charles E.\ Leiserson; 
Ronald L.\ Rivest; and Clifford Stein, {\it Introduction to Algorithms,} 
3rd edition, MIT Press, 2009. 


\bibitem[CKMW08]{krick1} Felipe Cucker, Teresa Krick, Gregorio Malajovich, 
and Mario Wschebor, {\it ``A numerical algorithm for zero counting. 
I: Complexity and Accuracy,''} J.\ Complexity, Vol.\ 24 (2008), pp.\ 
582--605. 

\bibitem[CKS18]{krick2} Felipe Cucker, Teresa Krick, Michael Shub, 
{\it ``Computing the homology of real projective sets,''} Found.\ Comput.\ 
Math.\ (2018) 18: 929-970. 

\bibitem[DGRM20]{sparsechem1} Alicia Dickenstein; Magal\'{\i} Giaroli; 
Rick Rischter; Mercedes P\'{e}rez Mill\'{a}n, {\it ``Parameter regions that 
give rise to $2[n/2]+1$ positive steady states in the $n$-site phosphorylation 
system,''} Mathematical Biosciences and Engineering, 2019, 16(6):7589--7615. 

\bibitem[DMST19]{sparsebio2} Alicia Dickenstein; Mercedes P\'{e}rez Mill\'{a}n; 
Anne Shiu; and Xiaoxian Tang, {\it ``Multistationarity in Structured Reaction 
Networks,''} Bulletin of Mathematical Biology (2019) 81(5), 1527-1581. 


\bibitem[DKdW18]{timo} Mareike Dressler, Adam Kurpisz, and Timo de Wolff, 
{\it ``Optimization over the boolean hypercube via sums of nonnnegative 
circuit polynomials,''} Math ArXiV preprint {\tt arXiv:1802.10004} . 

\bibitem[DS41]{ds} R.\ J.\ Duffin and A.\ C.\ Schaeffer, {\it ``A refinement 
of an inequality of the brothers Markoff,''} Transactions of the American 
Mathematical Society, {\bf 50}, pp.\ 517--528, 1941.  

\bibitem[Emi20]{emiris} Ioannis Z.\ Emiris, {\tt MixedVolume-SparseResultants} 
software package,\\  
{\tt https://github.com/iemiris/MixedVolume-SparseResultants } . 

\bibitem[EPR19]{epr1} Alperen A.\ Erg\"{u}r, Grigoris Paouris, and J.\ Maurice 
Rojas, {\em ``Probabilistic Condition Number Estimates for Real Polynomial 
Systems I: A Broader Family of Distributions,''}  
Foundations of Computational Mathematics, Feb.\ 2019, Vol.\ 19, No.\ 1, pp.\ 
131--157. 

\bibitem[EPR20]{epr2} Alperen A.\ Erg\"{u}r, Grigoris Paouris, and J.\ Maurice 
Rojas, {\em ``Probabilistic Condition Number Estimates for Real Polynomial 
Systems II: Structure and Smoothed Analysis,''} Math ArXiV preprint 
{\tt arXiv:1809.03626 } . 

\bibitem[For19]{forsgard} Jens Forsg\aa{}rd, {\it ``Defective dual varieties 
or real spectra,''} Journal of Algebraic Combinatorics, February 2019, Volume 
49, Issue 1, pp.\ 49--67. 

\bibitem[vzGG13]{vzg} Joachim von zur Gathen and J\"urgen Gerhard,
{\it Modern Computer Algebra,} 3rd ed., Cambridge University
Press, 2013.


\bibitem[GKZ94]{gkz94} Israel M.\ Gel'fand, Mikhail M.\ Kapranov, and
Andrei V.\ Zelevinsky, {\it Discriminants, Resultants and Multidimensional
Determinants,} Birkh\"auser, Boston, 1994. 

\bibitem[Gra99]{grabiner} David J.\ Grabiner, {\it ``Descartes' Rule of 
Signs: Another Construction,''} The American Mathematical Monthly, 
Vol.\ 106, No.\ 9 (Nov.\ 1999), pp.\ 854--856. 

\bibitem[Gre16]{grenet} Bruno Grenet, {\it ``Bounded-degree factors of 
lacunary multivariate polynomials,''} J.\ Symb.\ Comput.\ 75, pp.\ 171--192, 
2016 (ISSAC 2014 special issue).    

\bibitem[Har80]{hardt} Robert M.\ Hardt, {\it ``Semi-Algebraic Local-Triviality 
in Semi-Algebraic Mappings,''} American Journal of Mathematics, Vol.\ 102, 
No.\ 2 (Apr., 1980), pp.\ 291--302. 

\bibitem[HvdH19]{harvey} David Harvey and Joris van der Hoeven, {\em ``Integer 
multiplication in time $O(n\log n)$,''} 
HAL preprint {\em https://hal.archives-ouvertes.fr/hal-02070778 } . 

\bibitem[Her51]{hermite} Charles Hermite, {\it ``Sur l'introduction des 
variables continues dans la th\'{e}orie des nombres,''}  
J.\ Reine Angew.\ Math., 41:191--216, 1851. 


\bibitem[HS95]{hs} Birkett Huber and Bernd Sturmfels, 
{\it ``A Polyhedral Method for Solving Sparse Polynomial Systems,''} Math.\
Comp.\ {\bf 64} (1995), no.\ 212, pp.\ 1541--1555.


\bibitem[Kho80]{kho} Askold G.\ Khovanski{\u\i}, {\it ``A class of systems of 
transcendental equations,''} Dokl.\ Akad.\ Nauk SSSR 255 (1980), no.\ 4, 
pp.\ 804–807.

\bibitem[Kho91]{few} Askold G.\ Khovanski{\u\i}, {\it Fewnomials,}
AMS Press, Providence, Rhode Island, 1991.  

\bibitem[KPT15a]{lostfew} Pascal Koiran; Natacha Portier; and Sebastian 
Tavenas, {\it ``On the intersection of a sparse curve and a low-degree curve: 
A polynomial version of the lost theorem,''} Discrete and Computational 
Geometry, 53(1):48-63, 2015.  

\bibitem[KPT15b]{realtau} Pascal Koiran; Natacha Portier; and Sebastian 
Tavenas, {\it ``A Wronskian approach to the real tau-conjecture,''} 
Journal of Symbolic Computation, 68(2):195-214, 2015. 


\bibitem[Kro95]{kronecker} Leopold Kronecker, {\it Werke,} Vol.\ 1, Leipzig, 
Teubner (1895). 

\bibitem[Kus77]{kush} Anatoly Georgievich Kushnirenko, {\it ``Newton
Polytopes and the B\'ezout Theorem,"} Functional Analysis and its Applications 
(translated {}from Russian), vol.\ 10, no.\ 3, July--September (1977), 
pp.\ 233--235.

\bibitem[LL11]{leeli} Tsung-Lin Lee and Tien-Yien Li, {\it ``Mixed
volume computation in solving polynomial systems,''}
in Randomization, Relaxation,
and Complexity in Polynomial Equation Solving, Contemporary Mathematics,
vol.\ 556, pp.\ 97--112, AMS Press, 2011.

\bibitem[Leg14]{legall} Fran\c{c}ois Legall, {\it ``Powers of tensors and fast 
matrix multiplication,''} Proceedings of ISSAC ({I}nternational 
{S}ymposium on {S}ymbolic and {A}lgebraic {C}omputation) 2014, ACM Press, 
pp.\ 296--303, 2014. 




\bibitem[LRW03]{tri} Tien-Yien Li; J.\ Maurice Rojas; and 
Xiaoshen Wang, {\it ``Counting Real Connected Components of Trinomial 
Curves Intersections and m-nomial Hypersurfaces,''} Discrete and 
Computational Geometry, 30:379--414 (2003).    


\bibitem[Lio51]{liouville} Joseph Liouville, {\it ``Sur des classes tr\`{e}s 
\'{e}tendues de quantit\'{e}s dont la valeur n'est ni alg\'{e}brique, ni 
m\^{e}me r\'{e}ductible \'{a} des irrationnelles alg\'{e}briques,''}  
Journal Math.\ Pures et Appl., {\bf 16} (1851), pp.\ 133-142.  

\bibitem[Mah64]{mahler} Kurt Mahler, {\it ``An inequality for the discriminant 
of a polynomial,''} The Michigan Mathematical Journal, 11(3):257–262, 1964. 

\bibitem[Mar89]{markov} A.\ A.\ Markov, {\it ``On a certain problem of 
D.\ I.\ Mendeleiff,''} (in Russian) 
Utcheniya Zapiski Imperatorskoi Akademii Nauk, {\bf 
62}, pp.\ 1--24, 1889. 

\bibitem[Mas85]{masser} David W.\ Masser, {\it ``Open Problems,''}
Prod.\ Symp.\ Analytic Number Theory (ed.\ by
W.\ W.\ L.\ Chen), Imperial Coll.\ London, 1985.

\bibitem[Mat00]{matveev} E.\ M.\ Matveev, {\it ``An explicit lower bound for a 
homogeneous rational linear form in logarithms of algebraic numbers, II''}, 
Izv.\ Ross.\ Akad.\ Nauk Ser.\ Mat.\ {\bf 64} (2000), pp.\ 125--180; 
English transl.\ in Izv.l\ Math.\ {\bf 64} (2000), pp.\ 1217–1269.

\bibitem[Mig82]{mignotte} Maurice Mignotte, {\it ``Some Useful Bounds,''} 
Computing, Suppl.\ 4, pp.\ 259--263 (1982), Springer Verlag. 

\bibitem[Mon20]{mondal} Pinaki Mondal, {\it ``How many zeroes? Counting the number of solutions of systems of polynomials via geometry at infinity,''} 
Math ArXiV preprint {\tt 1806.05346}  

\bibitem[Nes03]{nesterenko} Yuri Nesterenko, {\it ``Linear forms in
logarithms of rational numbers,''} Diophantine approximation (Cetraro, 2000), 
pp.\ 53--106, Lecture Notes in Math., 1819, Springer, Berlin, 2003. 

\bibitem[Nit13]{nitaj} \scalebox{.95}[1]{Abderrahmane Nitaj, {\em The $abc$ 
Conjecture Home Page,} {\tt https://nitaj.users.lmno.cnrs.fr/abc.html}}  

\bibitem[Oes88]{oesterle} Joseph Oesterl\'e, {\it ``Nouvelles
approches du `Th\'eor\`eme' de Fermat,''} Ast\'erisque 161-2
(1988), pp.\ 165--186.

\bibitem[PPR19]{ppr} 
Grigoris Paouris, Kaitlyn Phillipson, and J.\ Maurice Rojas, 
{\em ``A Faster Solution to Smale's 17th Problem I: Real Binomial Systems,''}  
in proceedings of ISSAC 2019 (July 15-18, 2019, Beihang University, Beijing, 
China), ACM Press, 2019.

\bibitem[PRT09]{snc} Philippe P\'ebay; J.\ Maurice Rojas; and 
David C.\ Thompson, {\it ``Optimization and $\np_\R$-completeness of
certain fewnomials,''} proceedings of SNC 2009 (August 3--5, 2009, Kyoto,
Japan), pp.\ 133--142, ACM Press, 2009. 

\bibitem[PRS93]{proy} Paul Pedersen, Marie-Fran\c{c}oise Roy, and Aviva 
Szpirglas, {\it ``Counting real zeros in the multivariate case,''} in 
proceedings of Computational algebraic geometry (Nice, 1992), pp.\ 203–224, 
Progr.\ Math., 109, Birkh\"{a}user Boston, Boston, MA, 1993.

\bibitem[PR13]{pr} Kaitlyn Phillipson and J.\ Maurice Rojas, 
{\it ``Fewnomial Systems with Many Roots, and an Adelic Tau
Conjecture,''} in proceedings of Bellairs workshop on tropical and 
non-Archimedean geometry (May 6--13, 2011, Barbados), Contemporary 
Mathematics, vol.\ 605, pp.\ 45--71, AMS Press, 2013.  

\bibitem[Pra04]{prasolov} Victor V.\ Prasolov, {\it Problems and Theorems
in Linear Algebra,} translations of mathematical monographs, vol.\ 134,
AMS Press, 2004. 

\bibitem[RS02]{rs} Qazi Ibadur Rahman and Gerhard Schmeisser, {\it
Analytic Theory of Polynomials,} London Mathematical Society
Monographs 26, Oxford Science Publications, 2002.

\bibitem[Ren92]{renegar} Jim Renegar, {\it `` On the Computational Complexity
and Geometry of the First-Order Theory of the Reals, I--III,''}
J.\ Symbolic Comput.\ 13 (1992), no.\ 3, pp.\ 255--352. 

\bibitem[Roj99]{tgcp} J.\ Maurice Rojas, {\it ``Solving degenerate sparse 
polynomial systems faster,''} J.\ Symbolic Comput.\ 28 (1999), no. 1-2, 
pp.\ 155--186. 

\bibitem[Roj03]{why} J.\ Maurice Rojas, {\it ``Why Polyhedra Matter in
Non-Linear Equation Solving,''} Contemporary Mathematics, vol.\ 334, 
pp.\ 293--320, AMS Press, 2003.

\bibitem[Rou99]{rouillier} Fabrice Rouillier, {\it ``Solving zero-dimensional 
systems through the rational univariate representation,''} 
Appl.\ Algebra Engrg.\ Comm.\ Comput.\ 9 (1999), no.\ 5, pp.\ 433--461.  

\bibitem[Sag11]{sagraloff} Michael Sagraloff, {\it ``A General Approach to 
Isolating Roots of a Bit-stream Polynomial,''}  
Mathematics in Computer Science {\bf 4}, 481 (2010), 
Springer-Verlag. 

\bibitem[Sal76]{salamin} Eugene Salamin, {\it ``Computation of $\pi$
using arithmetic-geometric mean,''}  Math.\ Comput., 30 (1976), pp.\ 565--570

\bibitem[Sch86]{schrijver} Alexander Schrijver, {\it Theory of Linear
and Integer Programming,} John Wiley \& Sons, 1986. 

\bibitem[Sch80]{schwartz} Jacob T.\ Schwartz, {\it ``Fast Probabilistic
Algorithms for Verification of Polynomial Identities,''}
J.\ of the ACM 27, 701--717, 1980.

\bibitem[SL54]{descartes} David Eugene Smith and Marcia L.\ Latham, {\it
The Geometry of Ren\'e Descartes,} translated from the French and Latin
(with a facsimile of Descartes' 1637 French edition),
Dover Publications Inc., New York (1954).

\bibitem[Sto00]{storjophd} Arne Storjohann, {\it ``Algorithms for
Matrix Canonical Forms,''} doctoral dissertation, Swiss Federal
Institute of Technology, Zurich, 2000. 


\bibitem[V-W12]{virgi} Virginia Vassilevska-Williams, {\it ``Multiplying 
matrices faster than Coppersmith-Winograd,''} 
Proceedings of STOC (ACM Symposium on Theory of Computation) 2012, 
ACM Press, pp.\ 887--898, 2012.  

\bibitem[Ver10]{verschelde} Jan Verschelde, {\it ``Polynomial Homotopy 
Continuation with PHCpack'',} ACM Communications in Computer Algebra 
44(4):217-220, 2010.  

\end{thebibliography}

\end{document}